\documentclass[12pt]{amsart}
\usepackage{amssymb,latexsym,amscd,verbatim}
\usepackage[all]{xy}

\usepackage{mathrsfs}

\swapnumbers
\newtheorem{thm}[subsubsection]{Theorem}
\newtheorem{propose}[subsubsection]{Proposition}
\newtheorem{lemma}[subsubsection]{Lemma}

\newtheorem{cor}[subsubsection]{Corollary}
\theoremstyle{definition}
\newtheorem{defn}[subsubsection]{Definition}
\newtheorem{remark}[subsubsection]{Remark}

\newtheorem{conj}[subsubsection]{Conjecture}
\numberwithin{equation}{section}

\renewcommand{\d}{\mbox{\LARGE $\cdot $}}

\newcommand{\Gal}{{\rm Gal}\,}       
\newcommand{\Spec}{\operatorname{Spec}} 
\newcommand{\Hom}{\operatorname{Hom}}      
\newcommand{\Ext}{\operatorname{Ext}}      
\newcommand{\DM}{{\sf DM}}          
\newcommand{\M}{{\rm M}}   
\newcommand{\HI}{{\sf HI}} 
\newcommand{\PST}{{\rm PST}} 

\newcommand{\pu} {\underline{\pi}}

\newcommand{\ihom}{{\rm\underline{Hom}}}  

\newcommand{\Aff}{\mathbb{A}}   

\newcommand{\C}{\mathbb{C}}     
\newcommand{\Q}{\mathbb{Q}}     
\newcommand{\Z}{\mathbb{Z}}     
\newcommand{\N}{\mathbb{N}}
\newcommand{\G}{\mathbb{G}}     
\newcommand{\HH}{\mathbb{H}}    
\newcommand{\R}{\text{\rm R}}     
\renewcommand{\L}{\text{\rm L}}     

\newcommand{\im}{\operatorname{Im}}        
\renewcommand{\ker}{\operatorname{Ker}}  
\newcommand{\coker}{\operatorname{Coker}} 
\newcommand{\Pic}{{\rm Pic}}     
\newcommand{\Alb}{{\rm Alb}}     
\newcommand{\RPic}{{\rm RPic}}     
\newcommand{\LAlb}{{\rm LAlb}}     

\newcommand{\NS}  {{\rm NS}}      
\newcommand{\CH}{{\rm CH}}

\newcommand{\by}[1]{\stackrel{#1}{\rightarrow}}
\newcommand{\longby}[1]{\stackrel{#1}{\longrightarrow}}

\newcommand{\eproof}{\hfill$\Box$\\}
\renewcommand{\tilde}{\widetilde}
\newcommand{\df}{\mbox{\,${:=}$}\,}
\newcommand{\ie}{{\it i.e.\/},\ }
\newcommand{\cf}{{\it cf.\/}\ }
\newcommand{\eg}{{\it e.g.\/},\ }
\newcommand{\et} {\mbox{\scriptsize{\rm {\'e}t}}}

\newcommand{\eff}{\mbox{\scriptsize{\rm eff}}}
\newcommand{\gm}{\mbox{\scriptsize{\rm gm}}}

\newcommand{\into}{\hookrightarrow}

\newcommand{\colim}[1]{\mathop{\rm
Colim}_{\buildrel\over{#1}}}
\renewcommand{\lim}[1]{\mathop{\rm
Lim}_{\buildrel\over{#1}}}
\newcommand{\holim}[1]{\mathop{\rm
Holim}_{\buildrel\over{#1}}}

\newcommand{\boxtensor}{\def\boxtimesten{\Box\kern-7.59pt\raise1.2pt
\hbox{$\times$} }}                                  

\newcounter{elno}                   

\newcommand{\cA}{\mathcal{A}}
\newcommand{\cB}{\mathcal{B}}

\newcommand{\cE}{\mathcal{E}}
\newcommand{\cF}{\mathcal{F}}
\newcommand{\cG}{\mathcal{G}}

\newcommand{\cI}{\mathcal{I}}
\newcommand{\cJ}{\mathcal{J}}
\newcommand{\cK}{\mathcal{K}}
\newcommand{\cL}{\mathcal{L}}

\newcommand{\cN}{\mathcal{N}}

\newcommand{\cQ}{\mathcal{Q}}

\newcommand{\cT}{\mathcal{T}}

\newcommand{\catA}{\mathbf{A}}
\newcommand{\catB}{\mathbf{B}}
\newcommand{\catM}{\mathbf{M}}

\begin{document}

\title{1-motivic sheaves and the Albanese functor}
\author{Joseph Ayoub}
\address{Institut f\"ur Mathematik\\ Universit\"at Z\"urich\\Winterthurerstr. 190\\ CH-8057 Z\"urich\\ 
Switzerland}
\email{joseph.ayoub@math.uzh.ch}
\author{Luca Barbieri-Viale}
\address{Dipartimento di Matematica ``F. Enriques", Universit{\`a} degli Studi di Milano\\ Via C. Saldini, 50\\ I-20133 Milano\\ Italy}
\email{barbieri@math.unipd.it}
\begin{abstract} \begin{sloppypar} 
Using sheaf theoretic methods, we define functors 
$\L\pi_0 : \DM_{\eff}(k) \to D(\HI_{\leq 0}(k))$ and 
$\LAlb : \DM_{\eff}(k) \to D(\HI_{\leq 1}(k))$. The functor $\LAlb$
extends the one in \cite{BK} to non-necessarily geometric motives.
These functors are then used to define higher N\'eron-Severi groups and
higher Albanese sheaves. 
\end{sloppypar}
\end{abstract}
\maketitle

\tableofcontents

\section*{Introduction}

For a field $k$, say perfect, and a Grothendieck topology $\tau$ on the category of smooth $k$-schemes, such as the Nisnevich or the \'etale topology, 
we denote by $Shv^{\tau}_{tr}(k)$ the abelian category of $\tau$-sheaves with transfers on $Sm/k$.
Following Voevodsky, we consider $\DM^{\tau}_{\eff}(k)$ the full subcategory of the derived category $D(Shv^{\tau}_{tr}(k))$ whose objects are the $\Aff^1$-local complexes, \ie Voevodsky's (effective) motivic complexes. We refer to 
\cite[\S 3]{V} and \cite[Lect. 14]{VL} for an outline of this theory.

Attached to a smooth $k$-scheme $X$ we then get the representable $\tau$-sheaf $\Z_{tr}(X)\in Shv^{\tau}_{tr}(k)$ and the homological motive $\M (X)\in\DM^{\tau}_{\eff}(k)$ given by the $\Aff^1$-localization of $\Z_{tr}(X)$; recall that the $\Aff^1$-localization functor $D(Shv^{\tau}_{tr}(k)) \to\DM^{\tau}_{\eff}(k)$ is left adjoint to the obvious inclusion $\DM^{\tau}_{\eff}(k)\subset D(Shv^{\tau}_{tr}(k))$.
The smallest triangulated subcategory of $\DM_{\eff}^{\tau}(k)$ containing $\M(X)$ for $X\in Sm/k$ and stable by direct summands, is called the category of geometric (or constructible) motives and will be denoted by $\DM_{\eff,\gm}^{\tau}(k)$. When $\tau$ is the Nisnevich topology or the cohomological dimension of $k$ is finite, we obtain exactly the subcategory of compact objects.

Under some hypotheses (\eg the exponent characteristic of $k$ is inverted or $k$ is perfect and $\tau={\rm Nis}$), we know that the canonical $t$-structure on $D(Shv^{\tau}_{tr}(k))$ restricts 
to a $t$-structure on $\DM^{\tau}_{\eff}(k)$ whose heart is the abelian category $\HI_{tr}^{\tau}(k)\subset Shv^{\tau}_{tr}(k)$ of the homotopy invariant $\tau$-sheaves with transfers. This follows immediately from \cite[Th. 14.11]{VL}. For the \'etale topology, see \cite[D.3.3]{BK}. This $t$-structure is the so called homotopy $t$-structure.

\subsection{To the core} Notably, we may consider the triangulated subcategory $\DM^{\tau}_{\leq n}(k)\subset \DM^{\tau}_{\eff}(k)$ generated by  $\M (X)$ for $X$ of dimension $\leq n$ and closed with respect to direct sums, \ie the so called triangulated category of $n$-motivic complexes or $n$-motives.

A first step in the study of these subcategories was done by Voevodsky \cite[3.4]{V}: for example, one can see that the inclusion $\DM^{\rm Nis}_{\leq n}(k)\subset \DM_{\eff}^{\rm Nis}(k)$ has a right adjoint for all $n\geq 0$. Defining $\DM^{\tau}_{\leq n,\gm}(k)\subset \DM^{\tau}_{\leq n}(k)$ as before,
Voevodsky provided a description, rationally, of $\DM^{\rm Nis}_{\leq 0, \gm}(k)$ and $\DM^{\rm Nis}_{\leq 1, \gm}(k)$ in terms of Artin motives and Deligne $1$-motives (up to isogenies).

A second step was done by the second author jointly with B. Kahn, see \cite{BK}.
The category $\DM^{\et}_{\leq 1, \gm}(k)$ is described as the bounded derived category of Deligne 1-motives, for a suitable exact structure, after inverting  the exponential characteristic $p$ of the perfect field $k$,\ {\it via}\, a fully-faithful $\Z[1/p]$-linear embedding Tot into $\DM^{\et}_{\eff ,\gm}(k)$.
Furthermore, such embedding provides the homotopy $t$-structure on the derived category of Deligne $1$-motives whose heart is the $\Z[1/p]$-linear category of (constructible) $1$-motivic sheaves, see \cite[\S 3]{BK}.

A key result of \cite{BK} is that Tot has, rationally, a left adjoint which refines, integrally, to a functor $\LAlb $ on $\DM^{\rm Nis}_{\eff ,\gm}(k)$, the motivic Albanese triangulated functor. Dually, composing with (motivic) Cartier duality, one obtains the functor $\RPic$. Applied to the motive $\M(X)$ of an algebraic $k$-scheme $X$ these functors provide natural objects $\LAlb (X)$ and $\RPic (X)$ in $\DM^{\et}_{\leq 1, \gm}(k)$. An important application is in view of their $1$-motivic homology which is providing  the $1$-motives predicted by Deligne's conjecture. See the forthcoming second part of \cite{BK} for a proof of this conjecture (up to isogenies).

\subsection{Have a bird} The general goal of this paper is the study of the categories of $n$-motives by sheaf theoretic methods providing new algebraic invariants.\\

 In Section 1 we introduce the key notion of $n$-motivic $\tau$-sheaf, see \ref{n-motivic}. To do so, we first define (non-necessarily constructible) $n$-generated and strongly $n$-generated $\tau$-sheaves, see \ref{n-gener}. Roughly speaking, $n$-motivic $\tau$-sheaves are obtained from strongly $n$-generated $\tau$-sheaves by applying the functor 
 $h^{\tau}_0$ that takes a $\tau$-sheaf to a homotopy invariant one in a universal way. This functor is defined as the left adjoint of the inclusion  $\HI_{tr}^{\tau}(k)\subset Shv^{\tau}_{tr}(k)$.
An example of $n$-motivic $\tau$-sheaf is given by 
 $h_0^{\tau}(\Z_{tr}(X))$ with $X$ smooth of dimension less than $n$. 
We also show (under some mild hypotheses) that the category $\HI_{\leq n}^{\tau}(k)$ of 
$n$-motivic $\tau$-sheaves is a cocomplete abelian category, 
see \ref{HIn-abelian}.
For $n=1$ we show that this category is generated by lattices and semi-abelian group schemes. Actually, we show a structure theorem for $1$-motivic \'etale sheaves, see Theorem~\ref{struct-1-mot}, including finitely presented (or constructible)  $1$-motivic \'etale sheaves, see \ref{colim-fp-1-mot} and \cf \cite[\S 3.2]{BK}.

It is easy to see that $\HI_{\leq 0}^{\tau}(k)\cong Shv_{tr }^{\tau}(k_{\leq 0})$, the category of $\tau$-sheaves with transfers on $0$-dimensional smooth $k$-schemes. This yields a functor $$\pi_0: Shv_{tr }^{\tau}(k)\to \HI_{\leq 0}^{\tau}(k)$$ left adjoint to the inclusion $\HI_{\leq 0}^{\tau}(k)\subset Shv_{tr }^{\tau}(k)$, see \ref{omot}. 
With some more efforts, by taking a suitable colimit of Serre's Albanese schemes 
(\cf \cite{Serre2}), we obtain a functor
$$\Alb : Shv_{tr }^{\tau}(k)\to \HI_{\leq 1}^{\tau}(k)$$ left adjoint to the inclusion $\HI_{\leq 1}^{\tau}(k)\subset Shv_{tr }^{\tau}(k)$, see \ref{1mot-left-adj}. We denote by $(-)^{\leq n}$ the restriction of these functors to $ \HI_{tr}^{\tau}(k)$ and we conjecture that, at least rationally, the functors $$(-)^{\leq n}: \HI_{tr}^{\et}(k) \to  \HI_{\leq n}^{\et}(k) $$ exist also for $n\geq 2$.

We finally propose a conjectural framework (still for $n\geq 2$), remarkably linked to the Bloch-Beilinson conjectural filtration on zero-cycles, which permits a better understanding of the categories 
$\HI^{\et}_{\leq n}(k)$ and implies the existence of the functors $(-)^{\leq n}$ (see \ref{weak-BBF} and \ref{n-conjecture}).\\

 In Section 2 we construct functors $\L\pi_0$ and $\LAlb$ on $D(Shv^{\tau}_{tr}(k))$ as  ``true" derived functors of the functors $\pi_0$ and $\Alb$ defined in the previous section. In order to derive $\Alb$ we have to go through the proof that there are enough $\Alb$-admissible complexes, see \ref{existance-LF}. The key point here is that if $X$ is a smooth $k$-scheme which is affine and $\NS^1$-local, \ie the N\'eron-Severi geometrically vanishes, then $\Z_{tr}(X)$ is $\Alb$-admissible, see \ref{existence-enough-Alb-adm} for details. The so obtained $\LAlb$ factors through the $\Aff^1$-localization yielding a functor on $\DM^{\tau}_{\eff}(k)$: our main goal is then Theorem~\ref{thm-LAlb-main}.
As a by-product, we get, under some technical assumptions,
an equivalence of categories  
$\DM_{\leq n}^{\tau}(k)\simeq D(\HI_{\leq n}^{\tau}(k))$
for $n=0, \, 1$. See \ref{0-der-equi} and \ref{thm-LAlb-main} for a precise formulation.
 
 Note that $\L\pi_0$ and $\LAlb$ both take compact objects to compact objects so that $\LAlb$ is an extension of the one (in \cite{BK}) to non constructible motives.  We then show the non existence of left adjoints to $\DM^{\tau}_{\leq n}(k)\subset \DM^{\tau}_{\eff}(k)$ for $n\geq 2$ and set a conjecture linking $\HI_{\leq n}^{\tau}(k)$ to $\DM_{\leq n}^{\tau}(k)$.\\

In Section 3 we apply the functors $\L\pi_0$ and $\LAlb$ to the meaningful (non constructible) motivic complexes $\underline{\Hom}(\M(X), \Z (r)[2r])$ or $\underline{\Hom}(\Z(r)[2r],\M(X))$. The $s$-homology with respect to the homotopy $t$-structure is yielding a $0$-motivic sheaf whose group of $k$-points is the higher N\'eron-Severi $\NS^r(X, s)$. Similarly, we define the higher Picard $\Pic^r(X, s)$ and Albanese $\Alb_r(X,s)$ $1$-motivic sheaves.

The $\NS^r(X, 0)$ are related to higher codimension cycles in the following manner. Recall that the $h_0^{\rm Nis}$ of 
$\ihom(\M (X),\Z(r)[2r]))$ is the Nisnevich sheaf $\text{CH}_{/X}^r$ associated to the
presheaf $U\rightsquigarrow \text{CH}^r(U\times X)$ given by the Chow group of codimension $r$-cycles. Since $\pi_0(\text{CH}_{/X}^r) = \text{NS}_{/X}^r$ by the Theorem~\ref{equality-NS-pi-CH} we obtain that $\NS^r(X, 0)=\NS^r(X)$ is the classical N\'eron-Severi group of codimension $r$ cycles modulo algebraic equivalence.

\subsection*{\it Notation and Conventions} We let $k$ be our base field and $p$ its exponential characteristic. By scheme we always mean a finite type $k$-scheme. We warn the reader that all over in this paper we tacitly invert $p$ in the Hom groups of all categories constructed out of \'etale sheaves.

For the sake of exposition we here provide a comparison between some of the notations adopted in this paper and the corresponding existing notations in the book \cite{VL} as follows:

\medskip

\begin{center}
\begin{tabular}{|l|l|l|}
\hline
Paper & Book & Meaning\\
\hline
\hline
$Cor (k)$ & $ Cor_k$ & {\small category of finite correspondences}\\
\hline 
$Cor (X, Y)$ & $Cor (X, Y)$ & {\small group of finite correspondences from $X$ to $Y$} \\
\hline
$\PST (k)$ & $\mathbf{PST}(k)$ & {\small category of presheaves with transfers} \\
\hline
$Shv_{tr}^{\tau}(k)$ & $Sh_{\tau}(Cor_k)$ & {\small category of $\tau$-sheaves with transfers} \\
\hline
$\DM_{\eff}^{\tau}(k)$ & $\mathbf{DM}_{\tau}^{\eff}(k)$ & {\small Voevodsky category of  effective $\tau$-motives} \\
\hline
$\Z_{tr}(X) $ & $\Z_{tr}(X)$ & {\small representable presheaf with transfers} \\
\hline
\end{tabular}
\end{center}

\section{$n$-generated sheaves}
Let $Sm/k$ be the category of smooth schemes and $Cor(k)$ the category of finite correspondences of Voevodsky \cite[Lect. 1]{VL}.
Let $\tau\in\{ \text{co, Nis, {\'e}t}\}$ be one of the following Grothendieck topologies on $Sm/k$: coarse, Nisnevich or {\'e}tale topology.

\subsection{Generalities} 
Let $X\in Sm/k$. We denote by $\Z _{tr}(X)$ the representable  presheaf with transfers $$U \rightsquigarrow \Z _{tr}(X)(U) \df Cor (U, X)$$
For any presheaf with transfers $\cF$ we have by Yoneda:
\begin{equation}\label{reptr} \Hom (\Z _{tr}(X), \cF)= \cF (X)\end{equation}
Let $\PST (k)$ be the category of presheaves with transfers on $Sm/k$ and let $Shv_{tr}^{\tau}(k)$ be the full subcategory of $\tau$-sheaves. Recall that the presheaf $\Z _{tr}(X)$ is actually a $\tau$-sheaf (see \cite[Lemma 6.2]{VL}). 
Further denote by $\HI_{tr}^{\tau}(k)$ the full subcategory of homotopy invariant $\tau$-sheaves with transfers on $Sm/k$ (see \cite[Def. 2.15]{VL}).
 
\begin{lemma} \label{suslin}
The inclusions 
$$\HI_{tr}^{\tau}(k)\subset Shv_{tr}^{\tau}(k) \subset\PST (k)$$
admit left adjoints
$$ \PST (k)\longby{a_{\tau}} Shv_{tr}^{\tau}(k) \longby{ h_0^{\tau}}\HI_{tr}^{\tau}(k)$$
\end{lemma}

\begin{proof} 
The functor $a_{\tau}$ is the ``associated sheaf'' functor (\cf \cite{SGA4}).
Here we use that the $\tau$-sheaf $a_{\tau}(F)$ associated to a presheaf with transfers $F$ admits a unique structure of presheaves with transfers such that $F\to {\rm a}_{\tau}(F)$ is a morphism of presheaves with transfers. See \cite[3.1.4]{Voe1} for a proof in the case of Nisnevich topology and \cite[Th. 6.17]{VL} in the case of \'etale topology.

For $\tau = {\rm co}$ we let  $h^{\rm co}_0\df h_0$ the associated homotopy
 invariant presheaf functor, \ie the $H_0$ of the Suslin complex
 $C_*$. For $\tau \neq {\rm co}$ we define inductively (as in \cite{Deglise}), $h_0^{\tau, 0}\df {\rm id}$ and for all non negative integers
 $n\geq 0$ 
$$h^{\tau, n+1}_0 \df a_{\tau} h_0 h_0^{\tau, n}$$
and then take the colimit (in the category of presheaves)
$$h^{\tau}_0\df \colim{n\geq 0} h^{\tau, n}_0 $$

To show that $h_0^{\tau}$ takes values in the category of homotopy invariant $\tau$-sheaves, consider the following commutative diagram:
$$\xymatrix{h_0^{\tau,n}  \ar[r] \ar[d] & h_0^{\tau,n+1}  \ar[r] \ar[d] & h_0^{\tau,n+2}  \ar[d] \\ 
h_0h_0^{\tau,n}  \ar[r] \ar[d] & h_0 h_0^{\tau,n+1}  \ar[r] \ar[d] & h_0 h_0^{\tau,n+2}  \ar[d] \\ 
a_{\tau}h_0h_0^{\tau,n}  \ar[r] \ar[d] & a_{\tau} h_0 h_0^{\tau,n+1}  \ar[r] \ar[d] & a_{\tau} h_0 h_0^{\tau,n+2}  \ar[d] \\
h_0h_0^{\tau,n+1}  \ar[r]  & h_0h_0^{\tau,n+2}  \ar[r]  & h_0h_0^{\tau,n+3} }$$
Passing to the colimit we get the following sequence (using that $h_0$ commutes with colimits of presheaves):
$$\xymatrix{h_0^{\tau}  \ar[r] \ar@/^1pc/@{=}[rr]^{\;} & h_0 h_0^{\tau} \ar[r] \ar@/_1pc/@{=}[rr]_{\;}  & h_0^{\tau}  \ar[r] & h_0h_0^{\tau}}$$
which proves that $h_0^{\tau}=h_0h_0^{\tau}$. But $h_0^{\tau}(?)$ is a $\tau$-sheaf (because the topology $\tau$ is quasi-compact) and $h_0h_0^{\tau}(?)$ is homotopy invariant.

It is  easy to see that $h^{\tau}_0$ is a left adjoint, \eg note that on a homotopy invariant $\tau$-sheaf $\dag$ we get $h_0^{\tau}(\dag )= \dag$. 
\end{proof}

Under some mild hypotheses, we have $h_0^{\tau}=h_0^{\tau,1}=a_{\tau}h_0$ as the following proposition shows.

\begin{propose}\label{prop-abelian-hitr}
Assume that one of the following conditons:
\begin{enumerate}

\item $\tau={\rm co}$ is the coarse topology,

\item $k$ is perfect and $\tau={\rm Nis}$ is the Nisnevich topology,

\item the exponent characteristic $p$ of $k$ is inverted.

\end{enumerate}
Let $\cF$ be a homotopy invariant presheaf with transfers. Then $a_{\tau}(\cF)$ is strictly homotopy invariant, \ie ${\rm H}^n_{\tau}(-,a_{\tau}(\cF))$ is homotopy invariant for all $n$.
\end{propose}

\begin{proof}
When $\tau={\rm co}$ there is nothing to prove. 
For $\tau={\rm Nis}$ and $k$ perfect, this follows from \cite[Lect. 22]{VL}. 
If $k$ is not perfect, let $k_{insp}$ be the biggest totaly inseparable extension of $k$ (contained in an algebraic closure of $k$).
As remarked by Suslin, the base-change functor 
$Cor(k)\to Cor(k_{insp})$ becomes an equivalence of categories
when $p$ is inverted. It is then possible to extend 
Voevodsky's result to non perfect fields up to $p$-torsion.

Suppose that $\tau=\text{\'et}$ and $p$ is inverted. The following argument is similar to 
\cite [Lemma D.1.3]{BK}. 
Using the Hochschild-Serre spectral sequence, we may reduce to the case $k$ separably closed. 
Let $\cF_{tor}$ be the torsion sub-presheaf of $\cF$. By Suslin rigidity theorem 
\cite[Th. 7.20]{VL}, we know that $a_{\et}(\cF_{tor})$ is a constant \'etale sheaf (as $k$ is separably closed). By \cite[XV, Cor. 2.2]{SGA4}, we deduce that 
 $a_{\et}(\cF_{tor})$ is strictly homotopy invariant. Using the long exact sequence of cohomology, we reduce to the case of 
 $\cF'=\cF/\cF_{tor}$. 

Let $\cF''=\cF\otimes \Q/\cF'$. Using again Suslin rigidity theorem 
\cite[Th. 7.20]{VL} and the long exact sequence of cohomology we reduce to the case of $\cF\otimes \Q$. But if $\cG$ is a homotopy invariant presheaf with transfers taking values in the category of $\Q$-vector spaces, we have 
$a_{\et}(\cG)=a_{\rm Nis}(\cG)$ and ${\rm H}^n_{\et}(-,a_{\et}(\cG))
={\rm H}_{\rm Nis}^n(-,a_{\rm Nis}(\cG))$. The claim now follows from 
\cite[Lect. 22]{VL}.
\end{proof}

\begin{cor}\label{cor-ab-hitr}
Same assumption as in Proposition \ref{prop-abelian-hitr}.
The category $\HI_{tr}^{\tau}(k)$ is abelian complete and
cocomplete, the inclusion 
$\HI^{\tau}_{tr}(k)\subset Shv_{tr}^{\tau}(k)$ is exact and $h_0^{\tau}$ is right exact. 
\end{cor}

\subsubsection{} For $X\in Sm/k$ we let
$$h_0^{\tau}(X)\df h_0^{\tau}(\Z_{tr}(X))$$
For a homotopy invariant $\tau$-sheaf  $\cF\in \HI_{tr}^{\tau}(k)$ we thus obtain
\begin{equation}\label{rephi}
\Hom (h_0^{\tau}(X), \cF)= \Hom (\Z_{tr}(X), \cF) = \cF (X)
\end{equation}

\subsubsection{ } For $\cF\in  \PST (k)$ we  have a canonical map
$$\colim{X\to \cF} \Z _{tr} (X)\to \cF$$
where the colimit is taken over the category $Cor (k)/\cF$ whose objects are the elements in $\cF (X)$ for $X\in Sm/k$ or equivalently
(by \eqref{reptr}) maps of presheaves with transfers $\Z_{tr}(X) \to \cF$. 
Morphisms in $Cor(k)/\cF$ are commutative triangles of presheaves with transfers
$$\xymatrix{\Z_{tr}(X) \ar[r] \ar@/_/[dr] & \Z_{tr}(Y) \ar[d] \\
& \cF}$$
Note that the indexing category is \emph{pseudo-cofiltered} in the sense that any two objects are the target of two arrows having the same domain. Indeed for the two objects $\Z_{tr}(X)\to \cF$ and $\Z_{tr}(Y)\to \cF$ we can take $\Z_{tr}(X\coprod Y) \to \cF$. 

\begin{lemma}\label{lim} 
For $\cF\in  Shv_{tr}^{\tau}(k)$ we have an isomorphism
$$\colim{X\to \cF} \Z _{tr} (X)\xymatrix{\ar[r]^-{\sim} &} \cF$$
where the colimit is equally computed in $\PST(k)$ or in $Shv^{\tau}_{tr}(k)$.
\end{lemma}

\begin{proof} 
This is a well known fact. For any presheaf $\cF'\in \PST(k)$
consider the composition: 
$$\xymatrix{\Hom(\underset{}{}\cF,\cF')\ar[r] &  \Hom(\underset{X \to \cF}{\rm Colim} \; \Z_{tr}(X) , \cF')
\ar[r]^-{\sim} &    \underset{X\to \cF}{\rm Lim} \Hom(\Z_{tr}(X),\cF') \ar@{=}[d]\\
&  &  \underset{X \to \cF}{\rm Lim} \cF'(X)}$$ 
By Yoneda we need to prove that this is an isomorphism. 
Elements of $\lim{X\to \cF}\cF'(X)$ are families of  $\alpha'\in \cF'(X)$ indexed by $\alpha \in \cF(X)$ and satisfying the following compatibility with correspondences:
for any $\beta \in \cF(Y)$ and $\gamma\in Cor(X,Y)$ such that $\alpha=\gamma^*(\beta)$ we have $\alpha'=\gamma^*(\beta')$.  
In other terms, $\lim{X \to \cF} \cF'(X)$ is exactly the set of families of functions $(f_X:\cF(X) \to \cF'(X))_X$ compatible with the action of correspondences. To prove that such a family is a morphism of $\tau$-sheaves with transfers we still need to verify that $f_X$ are linear maps. This follows immediately from the diagram:
$$\xymatrix{\cF(X)\oplus \cF(X) \ar@{=}[r] \ar[d]_-{f_X\oplus f_X} & \cF(X\coprod X) \ar[r]^-{\gamma^*} \ar[d]^-{f_{X\coprod X}} & \cF(X)\ar[d]^-{f_X}\\
\cF'(X)\oplus \cF'(X) \ar@{=}[r] & \cF'(X\coprod X) \ar[r]^-{\gamma^*}  & \cF'(X) }$$ 
Where $\gamma$ is the sum of the two obvious inclusions $X \subset X\coprod X$.
\end{proof}

\begin{remark} 
The argument in the proof works for any site with finite coproducts
and a topology for which the family of morphisms $X_i \to \coprod_i X_i$ is a covering for any finite family $(X_i)_{i\in I}$. 
\end{remark}

\begin{cor}
For $\cF\in \HI_{tr}^{\tau}(k)$ we have
$$ \colim{X\to \cF} h_0^{\tau}(X)
\xymatrix{\ar[r]^{\sim} & } \cF$$
Here the colimit is computed in the category $Shv^{\tau}_{tr}(k)$. 
\end{cor}

\begin{proof} 
The map in Lemma \ref{lim} factors as follows
$$ \colim{X\to \cF} \Z _{tr} (X) \xymatrix{\ar@{->>}[r] &} \colim{X\to \cF}  h_0^{\tau}(X)\xymatrix{\ar[r] &} \cF$$
where the first map is surjective and the composition is an isomorphism.
\end{proof}

\subsubsection{} Let  $(Sm/k)_{\leq n}$ be the category of  smooth
schemes of dimension $\leq n$ with the topology $\tau$ (remark that
the dimension is stable under $\tau$-covers). Denote $\sigma_n : Sm/k \to
(Sm/k)_{\leq n}$ the continuous map of sites in the sense of \cite[Def.\/ 1.42 ]{MV}, given by the obvious inclusion  $
(Sm/k)_{\leq n}\subset Sm/k$.  Note that a priori $\sigma_n$ is not a morphism of sites \ie the pull-back functor is not exact.  

Consider the full subcategory $Cor(k_{\leq n})$ of $Cor(k)$ whose objects are the same of $(Sm/k)_{\leq n}$. We let $\PST (k_{\leq n})$ be the category of presheaves with transfers on $(Sm/k)_{\leq n}$: these are the additive contravariant functors from $Cor (k_{\leq n})$ to the category of abelian groups.

For $X\in (Sm/k)_{\leq n}$ we let  $\Z _{\leq n}(X)\in \PST (k_{\leq n})$ denote the presheaf with transfers $$\Z _{\leq n}(X)(U) \df Cor (U, X)$$ given by finite correspondences. For any presheaf with transfers $\cF\in \PST (k_{\leq n})$ we have
\begin{equation}\label{reptrn} \Hom (\Z _{\leq n}(X), \cF)= \cF (X)\end{equation}
Note that the presheaf $\Z _{\leq n}(X)$ is a $\tau$-sheaf.
Denote by $Shv_{tr}^{\tau}(k_{\leq n})$ the subcategory of $\tau$-sheaves in $\PST (k_{\leq n})$. The same proof as for Lemma \ref{lim} gives:

\begin{lemma}
\label{limn}
For $\cF\in  Shv_{tr}^{\tau}(k_{\leq n})$ we have that
$$\colim{(X\to \cF)_{\leq n}} \Z _{\leq n} (X)\xymatrix{\ar[r]^-{\sim} &}  \cF$$
where the colimit is taken over the category $Cor(k_{\leq n})/\cF$.
\end{lemma}

\subsubsection{} We have a restriction functor on $\tau$-sheaves with transfers
$$\sigma_{n*}: Shv_{tr}^{\tau}(k) \xymatrix{\ar[r] & } Shv_{tr}^{\tau}(k_{\leq n} )$$
which is clearly exact.
\begin{lemma}\label{upstar} 
 The functor $\sigma_{n*}: Shv_{tr}^{\tau}(k)\to Shv^{\tau}_{tr}(k_{\leq n})$ has a left adjoint  
$$\sigma_n^*: Shv_{tr }^{\tau}(k_{\leq n})\xymatrix{\ar[r] & } Shv_{tr}^{\tau}(k)$$ 
which is given by 
$$\sigma_n^*(\cF) \df \colim{(X\to \cF)_{\leq n}} \Z _{tr} (X)$$
for $\cF\in Shv_{tr }^{\tau}(k_{\leq n})$. Here the colimit is
computed in $Shv_{tr}^{\tau}(k)$.
\end{lemma}

\begin{proof} 
In fact, for $\cF\in Shv_{tr }^{\tau}(k_{\leq n})$ and $\cF'\in Shv_{tr}^{\tau}(k)$ we have, by Lemma \ref{limn},
$$\Hom (\cF, \sigma_{n*}(\cF'))=  
\Hom \left(\colim{(X\to \cF)_{\leq n}} \Z _{\leq n} (X), \sigma_{n*}(\cF')\right)$$ 
which is 
$$\displaystyle\lim{(X\to \cF)_{\leq n}} \Hom(\Z_{\leq n}(X), \sigma_{n*}(\cF' ))$$
Since we clearly have $\Hom(\Z_{\leq n}(X), \sigma_{n*}(\cF' )) =\Hom(\Z_{tr}(X), \cF' ) = \cF' (X)$, for all $X\in (Sm/k)_{\leq n}$, \cf \eqref{reptrn}, we obtain:
$$\Hom (\cF, \sigma_{n*}(\cF')) = \Hom \left(\colim{(X\to \cF)_{\leq n}} \Z _{tr} (X), \cF'\right) =  \Hom (\sigma_n^*(\cF), \cF') $$
\end{proof}

\begin{defn} \label{n-gener}
A $\tau$-sheaf  $\cF\in Shv_{tr}^{\tau}(k)$ is 
 {\it $n$-generated} if the counit
 $$\sigma_n^*\sigma_{n*}(\cF) \xymatrix{\ar[r] & } \cF$$ is a surjection. When it is an isomorphism we say that $\cF$ is {\it strongly} $n$-generated.
We denote by $Shv_{\leq n}^{\tau}(k)$ the subcategory of strongly $n$-generated $\tau$-sheaves.
 \end{defn}

\begin{remark}\label{n-gen-and-topology}
The property of being (strongly) $n$-generated is compatible with the change of
topology. For example if $\cF$ is an $n$-generated  Nisnevich sheaf
then $a_{\et}\cF$ is an $n$-generated {\'e}tale sheaf. Indeed, we have
 $ a_{\et}\sigma_{n*}\simeq \sigma_{n*} a_{\et} $ and $\sigma_{n}^*a_{\et}\simeq 
 a_{\et} \sigma_n^*$. Beware that in the last formula, the first $\sigma_n^*$ stands for the inverse image on \'etale sheaves whereas the second one stands for the inverse image on Nisnevich sheaves.
\end{remark}

\begin{lemma}
\label{ngen-coker-ext}
The property of being (strongly) $n$-generated is stable by cokernels and extensions in the category of $\tau$-sheaves.
\end{lemma}

\begin{proof}
We do this only for extensions in the case of $n$-generated sheaves;
the other cases are simpler. The result follows from:
$$\xymatrix{& \sigma_n^* \sigma_{n*}  \cE \ar[r] \ar@{->>}[d] & \sigma_n^*\sigma_{n*} \cF \ar[r] \ar[d] & \sigma_n^* \sigma_{n*} \cG \ar[r] \ar@{->>}[d] & 0 \\
0 \ar[r] & \cE \ar[r] & \cF \ar[r] & \cG \ar[r] & 0  }$$
and a diagram chase.
\end{proof}

\begin{lemma} \label{unitn} 
The unit ${\rm id} \xymatrix{\ar[r]^-{\sim} &} \sigma_{n*}\sigma_n^* $
  is invertible.
\end{lemma}

\begin{proof} 
For $\Z_{\leq n}(X)$ and $X\in (Sm/k)_{\leq n}$, we have
 $\sigma_n^* \Z_{\leq n}(X) = \Z_{tr}(X)$ and $\sigma_{n*} \Z_{tr}(X) =  \Z_{\leq n}(X)$. It follows that $\Z_{\leq n}(X) \simeq \sigma_{n*}\sigma_n^*\Z_{\leq n}(X)$. Using Lemma \ref{limn}, we only need to show that $\sigma_n^*$ and $\sigma_{n*}$ commute with colimit.
 This is clear for $\sigma_n^*$ as it is a left adjoint. For $\sigma_{n*}$, we use that it commutes with colimits of presheaves and with sheafification.
\end{proof}

Note the following useful corollary:

\begin{cor}
\label{tricky-ngen}
Let $\cF$ be a $\tau$-sheaf with transfers on $Sm/k$. Denote by $\cN$ the kernel of $\xymatrix{\sigma_n^*\sigma_{n*} (\cF) \ar[r] & \cF}$. If $\cN$ is $n$-generated then it is zero. 
\end{cor}

\begin{proof}
As $\sigma_{n*}$ is exact, we have a left exact sequence:
$$\xymatrix{0 \ar[r] & \sigma_{n*}(\cN) \ar[r] & \sigma_{n*}\sigma_n^* \sigma_{n*} (\cF) \ar[r] & \sigma_{n*}(\cF)}$$
Using \ref{unitn} and that the composition:
$$\xymatrix{\sigma_n^* \ar[r]^-{\sim} &  \sigma_n^*\sigma_{n*}\sigma_n^*\ar[r] & \sigma_n^*}$$
 is the identity, we see that $\sigma_{n*}(\cN)=0$. But as $\cN$ is $n$-generated, we have a surjection: $\xymatrix{0=\sigma_n^* \sigma_{n*}(\cN) \ar@{->>}[r] & \cN}$.\end{proof}

\begin{propose}\label{image}
The functor $\sigma_n^*$ in Lemma \ref{upstar} takes values in the category $Shv_{\leq n}^{\tau}(k)$ and it induces an equivalence between $Shv_{tr }^{\tau}(k_{\leq n})$  and the category of  strongly $n$-generated sheaves. 
\end{propose}

\begin{proof} Everything follows from Lemma \ref{unitn}.  The essential image of  $\sigma_n^*$ consists of  strongly $n$-generated sheaves because we always have that the composition of
$$\xymatrix{\sigma_n^*\ar[r]^-{\sim} & \sigma_n^*\sigma_{n*}\sigma_n^*\ar[r] & \sigma_n^*}$$
 is the identity and we have that the first map is an isomorphism.
\end{proof}

\begin{remark} 
An example of strongly $n$-generated sheaf is $\Z_{tr} (X)$ for $X\in (Sm/k)_{\leq n}$. 
It follows that $h_0^{\tau}(X)$  is $n$-generated. However we don't expect this sheaf to be strongly $n$-generated for $n\geq 1$. 
We leave it as an open (possibly hard) problem to prove (or disprove)
that $h_0(C)$ is not strongly $1$-generated for an elliptic curve or
even for $\mathbb{G}_m$.   
\end{remark}

\begin{defn} \label{n-motivic}
 A homotopy invariant $\tau$-sheaf  $\cF\in \HI_{tr}^{\tau}(k)$ is {\it $n$-motivic}  if  
$$\xymatrix{h_0^{\tau}(\sigma_n^*\sigma_{n*}(\cF) )\ar[r]^-{\sim} &  h_0^{\tau}(\cF)=\cF}$$
is an isomorphism. We let $\HI_{\leq n}^{\tau}(k)$ be the full subcategory of homotopy invariant $n$-motivic $\tau$-sheaves. 
\end{defn}

\begin{remark}\label{exampleh}
By definition any $n$-motivic $\tau$-sheaf is the $h_0^{\tau}$ of a strongly $n$-generated $\tau$-sheaf. Conversely, if a $\tau$-sheaf $\cF$ is strongly $n$-generated then $h_0^{\tau}(\cF)$ is $n$-motivic. Indeed, we have the following commutative square of epimorphisms:
$$\xymatrix{\sigma_n^*\sigma_{n*} (\cF) \ar@{->>}[r] \ar[d]_-{\sim} & \sigma_n^*\sigma_{n*} h_0^{\tau}(\cF) \ar@{->>}[d]\\
\cF \ar@{->>}[r] & h_0^{\tau}(\cF)}$$
Applying $h_0^{\tau}$ we get:
$$\xymatrix{h_0^{\tau}\sigma_n^*\sigma_{n*} (\cF )\ar@{->>}[r] \ar[d]_-{\sim} \ar@/_/[dr]^-{\sim} & h_0^{\tau}\sigma_n^*\sigma_{n*} h_0^{\tau}(\cF) \ar@{->>}[d]\\
h_0^{\tau}(\cF) \ar@{=}[r] & h_0^{\tau}(\cF)}$$
Which proves that the arrow $\xymatrix{h_0^{\tau}\sigma_n^*\sigma_{n*} h_0^{\tau}(\cF) \ar[r]^-{\sim} & h_0^{\tau}(\cF) }$ is invertible. In particular the $\tau$-sheaves $h_0^{\tau}(X)$ are $n$-motivic for smooth $k$-varieties of dimension $\leq n$.
\end{remark}

\begin{lemma}
\label{nmot-coker-ext}
Same assumption as Proposition \ref{prop-abelian-hitr}. 
The property of being $n$-motivic is stable by cokernels and extensions in $\HI^{\tau}_{tr}(k)$.
\end{lemma}

\begin{proof}
Recall (Corollary \ref{cor-ab-hitr}) that $h_0^{\tau}$ is right exact being the left adjoint of an exact functor.
Then use the same diagram chase as in the proof of Lemma \ref{ngen-coker-ext} adding $h_0^{\tau}$ on the top line.  
\end{proof}

Denote by $inc: \HI^{\tau}_{tr}(k) \subset Shv_{tr}^{\tau}(k)$ the obvious inclusion. We have the following weaker version of Lemma \ref{unitn}:

\begin{lemma}
\label{unitnmot}
The two natural transformations:
$$\xymatrix{(\sigma_{n*} inc) \ar[r]^-{\sim} &  (\sigma_{n*} inc)(h_0^{\tau} \sigma_n^*)(\sigma_{n*} inc) \ar[r]^-{\sim} & 
(\sigma_{n*} inc)}$$
are invertible.

\end{lemma}

\begin{proof}
 As the composition of the two arrows of the lemma is the identity, we need only to show that the left hand side is surjective when applied to any $\cF \in \HI^{\tau}_{tr}(k)$. This follows from the commutative diagram:
$$\xymatrix{\sigma_{n*} inc (\cF) \ar[r]^-{\sim} \ar@/_/[dr] & \sigma_{n*} \sigma_n^*\sigma_{n*} inc (\cF) \ar@{->>}[d]\\
& (\sigma_{n*} inc) (h_0^{\tau}  \sigma_{n}^*) (\sigma_{n*} inc) (\cF)}$$
and Lemma \ref{unitn}.
\end{proof}

\begin{cor}
\label{HIn-abelian}
Same assumption as Proposition \ref{prop-abelian-hitr}. 
The category $\HI^{\tau}_{\leq n}(k)$ is abelian and cocomplete. The inclusion $\HI_{\leq n}^{\tau}(k)\subset \HI_{tr}^{\tau}(k)$ is right exact. 
\end{cor}

\begin{proof}
Let $f: \cF \to \cF'$ be a morphism between two $n$-motivic
sheaves. By Lemma \ref{nmot-coker-ext}, $\coker(f)$ is $n$-motivic so that
$\HI_{\leq n}^{\tau}(k)$ admits cokernels. 
The category $\HI_{\leq n}^{\tau}(k)$ admits also kernels that are
given by $h^{\tau}_0 \sigma_n^*\sigma_{n *} \ker(f)$. One easily
checks that the image and coimage agree by applying the conservative
(on $\HI_{\leq n}^{\tau}(k)$) functor $\sigma_{n*}$.
\end{proof}

\begin{remark}
For $\tau=\text{{\'e}t}$ and $p$ inverted, we believe that the inclusion $\HI_{\leq n}^{\et}(k)\subset \HI_{tr}^{\et}(k)$ is also left exact. However, this seems a difficult problem. See Corollary \ref{nmot-thick} for a conjectural proof relying on \ref{weak-BBF}.
\end{remark}

The following is a homotopy invariant version of Corollary \ref{tricky-ngen}:

\begin{cor}
\label{tricky-nmot}
Let $\cF$ be a homotopy invariant $\tau$-sheaf with transfers on $Sm/k$. Denote by $\cN$ the kernel of $\xymatrix{h_0^{\tau}\sigma_n^*\sigma_{n*} (\cF) \ar[r] & \cF}$. If $\cN$ is $n$-generated then it is zero. 
\end{cor}

\begin{proof}
As $\sigma_{n*}$ is exact, we have a left exact sequence:
$$\xymatrix{0 \ar[r] & \sigma_{n*}(\cN)\ar[r] & \sigma_{n*}h_0^{\tau}\sigma_n^* \sigma_{n*} (\cF )\ar[r] & \sigma_{n*}(\cF)}$$
By \ref{unitnmot}, we get $\sigma_{n*}(\cN)=0$. But as $\cN$ is $n$-generated, we have a surjection: $\xymatrix{0=\sigma_n^* \sigma_{n*}(\cN) \ar@{->>}[r] & \cN}$.
\end{proof}

\subsection{$0$-generated} Recall that a \emph{lattice} is a presheaf
which is representable by  a $k$-group scheme locally constant for the
{\'e}tale topology with geometric fiber isomorphic to a free finitely
generated abelian group. This is an example of $0$-generated {\'e}tale
sheaf.

\subsubsection{} 
For a reduced $k$-scheme $X$ one has the Stein factorization
$$\xymatrix{X \ar[r] & \pu_0(X) \ar[r] & \Spec(k)}$$
where $\pu_0(X)$ is the spectrum of the integral closure of $k$ in
$\Gamma(X,\mathcal{O}_X)$. If $X$ is smooth and $l$ is a finite \'etale extension of $k$, we have a canonical isomorphism
$$Cor(X,\Spec(l)) \simeq Cor(\pu_0(X),\Spec(l))\simeq \Z^{\pi_0(|X\otimes_k l|)}$$
where $|X\otimes_k l|$ is the Zariski topological space underlying the scheme $X\otimes_k l$ and $\pi_0(|X\otimes_k l|)$ is the set of connected components.

We thus have a functor
$\pu_0:\xymatrix{Cor (k) \ar[r] & Cor (k_{\leq 0})}$ which is left
adjoint to the inclusion $\sigma_0:Cor(k_{\leq 0}) \subset Cor(k)
$. The functor $\pu_0$ clearly induces a map of $\tau$-sites, so that we have a pair $(\pu_0^*, \pu_{0*})$ of adjoint functors:
$$\xymatrix{Shv_{tr}^{\tau}(k)   \ar[r]^{\pu_0^*}     &   \ar@/^1.7pc/[l]^{\pu_{0*}}  Shv_{tr}^{\tau}(k_{\leq 0}) }$$ From the adjunction $(\pu_0,\sigma_0)$, one  immediately gets an
adjunction $(\pu_0^*,\sigma_0^*)$. This gives a canonical isomorphism
$\pu_{0*}\simeq \sigma_0^*$.

\begin{lemma}\label{stronghom}
A strongly $0$-generated $\tau$-sheaf is homotopy
invariant. Furthermore, the functor
$\sigma_0^*:Shv_{tr}^{\tau}(k_{\leq 0}) \to Shv_{tr}^{\tau}(k)$
induces an equivalence of categories between
$Shv_{tr}^{\tau}(k_{\leq 0})$ and $\HI_{\leq 0}^{\tau}(k)$.
\end{lemma}

\begin{proof}
Take $\cF=\sigma_0^*\cF_0=\pu_{0*}\cF_0$. Using
$\cF(X)=\cF_0(\pu_0(X))$ we only need to show that $\pu_0(X\times_k
\Aff^1_k)=\pu_0(X)$ which is true, more generally, for $X$ reduced. The last assertion follows from Proposition \ref{image}.
\end{proof}

\begin{defn}
A $0$-motivic $\tau$-sheaf $\cE$ is \emph{finitely generated} if there exists an \'etale $k$-algebra $l$ and a surjection $\xymatrix{\Z_{tr}(\Spec(l))\ar@{->>}[r] & \cE}$.
\end{defn}

\begin{cor}
A $0$-motivic $\tau$-sheaf is a filtered colimit of finitely generated
$0$-motivic sheaves.  
\end{cor} 

\begin{proof}
Let us say that a $\tau$-sheaf with transfers $\cE_0$ on $(Sm/k)_{\leq 0}$ is finitely generated if there exist an \'etale $k$-algebra $l$ and a surjection $\xymatrix{\Z_{\leq 0}(\Spec(l))\ar@{->>}[r] & \cE_0}$.
By Lemma \ref{stronghom}, a $0$-motivic $\tau$-sheaf $\cE=\sigma_0^*\cE_0$ is finitely generated if and only if $\cE_0$ is finitely generated. We are thus reduced to prove the corresponding statement for $\tau$-sheaves with transfers on $(Sm/k)_{\leq 0}$.
But it is clear that such a $\tau$-sheaf $\cF_0$ is a filtered union of images of $\Z_{\leq 0}(\Spec(l))\to \cF_0$ with $l$ an \'etale $k$-algebra.
\end{proof}

\begin{cor} \label{PI-zero}
The embedding $\HI_{\leq 0} ^{\tau}(k)\into
  Shv_{tr}^{\tau}(k)$ has a  left adjoint 
$$\pi_0: Shv_{tr}^{\tau}(k)\to\HI^{\tau}_{\leq 0}(k)$$ 
given by
$$\pi_{0}(\cF)\df \colim{X\to \cF}\Z_{tr}(\pu_0(X))$$
\end{cor}

\begin{proof}
Indeed, $\HI_{\leq 0}^{\tau}(k)$ is simply the
subcategory of strongly $0$-generated $\tau$-sheaves which in turn is
equivalent to $Shv_{tr}^{\tau}(k_{\leq 0})$. 
Under this equivalence the inclusion $\HI_{\leq 0}^{\tau}(k)\subset Shv_{tr}^{\tau}(k)$ 
is given by $\sigma_0^*\simeq \pu_{0*}$. The latter admits $\pu_0^*$ as a left adjoint.
The formula follows from Lemma \ref{lim} and the commutation of left adjoints with colimits.
\end{proof}

\begin{defn} \label{omot} Denote by $(-)^{\leq 0}:\HI_{tr}^{\tau}(k) \to \HI_{\leq 0}^{\tau}(k)$ the restriction of $\pi_0$ to $\HI_{tr}^{\tau}(k) \subset Shv_{tr}^{\tau}(k)$. It is clearly the left adjoint of the inclusion $\HI_{\leq 0}^{\tau}(k) \subset \HI_{tr}^{\tau}(k) $.
\end{defn}

\begin{propose} \label{ogen}
Assume one of these conditions is fulfilled:
\begin{enumerate}

\item[(a)] $k$ is separably closed,

\item[(b)] $\tau$ is the {\'e}tale topology,

\item[(c)] that we work with rational coefficients. 
 \end{enumerate}
Then a $0$-generated $\tau$-sheaf is strongly $0$-generated and hence
$0$-motivic. The category $\HI_{\leq 0}^{\tau}(k)\subset
Shv_{tr}^{\tau}(k)$ is a Serre or thick abelian subcategory, \ie stable under
extensions, subobjects and quotients.
\end{propose}

\begin{remark}
\label{point-k-separ}
If $k$ is separably closed then any smooth $k$-scheme has a rational point.
\end{remark}

We first prove the following lemma:

\begin{lemma}
\label{boh}
Let $\cF$ be a $\tau$-sheaf.
Under one of the assumptions in Proposition \ref{ogen} the
morphism 
$$\xymatrix{\mathcal{F} \ar[r] & \pu_{0*} \pu_0^*\mathcal{F} }$$
is surjective.
\end{lemma}

\begin{proof} 
Using Lemma \ref{lim} we are left to show the statement for 
$$\xymatrix{\Z_{tr}(X) \ar[r] &  \Z_{tr}(\pu_0(X))}$$
This is clear when assuming (b) or (c). For (a), one uses Remark \ref{point-k-separ}.
\end{proof}

\subsubsection{} Let $\cF=\sigma_0^* \cF_0$ be a strongly 0-generated $\tau$-sheaf and suppose given a morphism $i:\xymatrix{\mathcal{E} \ar[r] & \cF}$. Because $\sigma_0^*\simeq \pu_{0*}$, this is equivalent to give a morphism $\xymatrix{\pu_0^*\mathcal{E} \ar[r] & \cF_0}$. We have by this a factorization:
$$\xymatrix{\mathcal{E} \ar[r] & \pu_{0*} \pu_0^*\mathcal{E} \ar[r] & \mathcal{F} }$$
By Lemma \ref{boh}, the first arrow is surjective. It follows that if $i$ is injective, we have an isomorphism $\mathcal{E} \simeq \pu_{0*}\pu_0^*\mathcal{E}\simeq \sigma_0^* \pu_0^*\mathcal{E}$. We have proven:

\begin{lemma}
\label{bohh}
Same assumption as in Proposition \ref{ogen}. Any subsheaf of a strongly $0$-generated $\tau$-sheaf is again strongly $0$-generated.
\end{lemma}

\smallskip

\noindent{\it Proof of Proposition~\ref{ogen}.}  
Let $\cF$ be a $0$-generated $\tau$-sheaf. 
By Proposition \ref{image}, $\sigma_0^*\sigma_{0*}(\cF)$ is strongly $0$-generated.
The kernel $\cN$ of the surjective morphism
$\xymatrix{\sigma_0^* \sigma_{0*}(\cF)\ar@{->>}[r] & \cF}$
is then $0$-generated by Lemma \ref{bohh}.
By Corollary \ref{tricky-ngen}, $\cN$ is zero.
The other claims are already proven in Lemma \ref{stronghom} and  Lemma \ref{ngen-coker-ext}.
\eproof

\begin{remark}
The category $\HI_{\leq 0}^{\et}$ is the smallest cocomplete Serre abelian subcategory of $Shv_{tr}^{\et}(k)$ containing lattices. Indeed $\HI_{\leq 0}^{\tau}(k)$ is equivalent to $Shv^{\tau}_{tr}(k_{\leq 0})$.  
\end{remark}

\begin{remark}
Under the assumption of Proposition \ref{ogen}, a subsheaf of a finitely generated $0$-motivic $\tau$-sheaf is again a finitely generated $0$-motivic $\tau$-sheaf as one easily checks by reducing to the case of $\tau$-sheaves with transfers on $(Sm/k)_{\leq 0}$. In particular, any finitely generated $0$-motivic $\tau$-sheaf $\cF$ admits a presentation: 
$$\Z_{tr}(\Spec(l_1))\to \Z_{tr}(\Spec(l_0))\to \cF\to 0$$
where $l_0$ and $l_1$
are \'etale $k$-algebras. Thus it makes sense to say that $\cF$ is finitely presented.
\end{remark}

\subsection{$1$-generated}  Let $G$ be a commutative  group scheme whose connected component of the identity $G^0$ is a semi-abelian variety and $\pi_0 (G)$ is finitely generated, \ie a semi-abelian scheme with torsion in the terminology of \cite[Def. 3.6.4]{BK}. Recall that a semi-abelian variety is an extension of an abelian variety by a torus. 
In the following we refer to such a $G$ as a semi-abelian group scheme for short.

Notably $G$ is a quotient of the Serre-Albanese scheme $\Alb(C)$ of a suitable smooth subvariety $C$ of $G$ of dimension $\leq 1$ (up to $p$-torsion). It follows that $G$ is $1$-generated as $\Alb(C)$ represents $h_0^{\et}(C)$  by Voevodsky \cite[Sect. 3.4]{V}.   

\subsubsection{Warnings and abuse of notation} From now on we stick to the case $\tau=\text{\'et}$ and invert the exponential characteristic $p$ of $k$.  All statements of this section hold only after inverting $p$. We will make the following abuse of notation: writing $\Z_{tr}(X)$ we mean $\Z[1/p]_{tr}(X)$ and writing $G$ we mean $G[1/p]$ in the corresponding $\Z[1/p]$-linear categories.

Note that given a smooth commutative group scheme $G$,  the \'etale sheaf $\underline{G}$ on $Sm/k$ represented by $G$ has a canonical structure of presheaf with transfers (\cf \cite[Lemma 1.3.2]{BK} and \cite{Samu}). This gives a functor from the category of smooth group schemes to the category of presheaves with transfers. One can easily prove that this functor is fully faithful. For this reason, we identify a smooth group scheme with the presheaf with transfers that represents. 

Further,  for an arbitrary sheaf $\cF\in Shv_{tr}^{\et}(k)$, we will denote $$\cF^0\df\ker(\cF \to \pi_0(\cF))$$
by making use of Corollary \ref{PI-zero}. We then say that such a sheaf $\cF$ is connected if $\pi_0(\cF)=0$.

\subsubsection{Serre-Albanese scheme} Recall by \cite{Serre2, Rama} that for a smooth $k$-variety $X$ we have a universal morphism $X \to \Alb(X)$ with $\Alb(X)$ a semi-abelian scheme as above. The group scheme $\Alb(X)$ is the Serre-Albanese scheme of $X$. The map $X \to \Alb(X)$ can be extended to a morphism of presheaves with transfers $\Z_{tr}(X) \to \Alb(X)$ (see \cite{BK, Samu}). 

\begin{lemma}\label{lemm-surj-to-alb}
Let $X$ be a smooth $k$-scheme. The morphism 
$\Z_{tr}(X)\to \Alb(X)$ is surjective for the \'etale topology.
\end{lemma}

\begin{proof}
Indeed, the image of $\theta:\Z_{tr}(X)\to \Alb(X)$ is an \'etale subsheaf of $\Alb(X)$. In particular, it is a homotopy invariant Nisnevich sheaf. To check that ${\rm Im}(\theta)=\Alb(X)$, we only need to look on function fields of smooth $k$-varieties (as follows from \cite[Lemma 22.8]{VL}). As ${\rm Im}(\theta)$ and $\Alb(X)$ are both \'etale sheaves, we may replace this function field by finite \'etale extensions. We are then reduced to show that $\Z_{tr}(X\otimes_k K)(K) \to \Alb(X\otimes_k K)(K)$ 
 is surjective for all extensions $K$ of $k$ that are separably closed. As, we invert the exponential characteristic of $k$, we may even suppose (using a transfers argument) that $K$ is algebraically closed. We are then left to show that the group of points of $\Alb(X)$ over an algebraically closed field is generated by the classes of closed points of $X$, which is a well-known fact.
\end{proof}

We now want to understand the subsheaves of $1$-motivic sheaves. Unfortunately, we can't use here the formalism of adjoint functors as in the previous paragraph; we are forced to give a direct proof of:

\begin{lemma}
\label{sub1mot}
Let $\cF$ be a $1$-motivic sheaf. Any subsheaf of $\cF$ is again $1$-motivic.
\end{lemma} 

\begin{proof} 
We break the proof in three steps. In the first two steps we show that a subsheaf of $\cF$ is $1$-generated. In the third part we deduce that this subsheaf is $1$-motivic.

\smallskip

\noindent
\emph{Step 1:} Consider first the case of $\cF=h_0^{\et}(C)$ with $C$ a smooth
scheme of dimension $\leq 1$, which is a $1$-motivic sheaf by  Remark \ref{exampleh}. 
Fix a subsheaf $\cE \subset \cF$. 
We can see $\cE$ as a filtered union of images of $\xymatrix{\Z_{tr}(X) \ar[r] & h_0^{\et}(C)}$. 
Actually, we may suppose that $\cE$ is the image of a map, \ie 
 $\cE=\text{Im}(a:\Z_{tr}(X) \to h_0^{\et}(C))$,
because any subsheaf is a filtered union of such images  and a colimit of $1$-generated $\tau$-sheaves is also $1$-generated. 

Since $h_0^{\et}(C)$ is represented by a semi-abelian group scheme $G$ then $a$ factors through $\Alb(X)$: 
$$\xymatrix{\Z_{tr}(X) \ar[r]^-{\theta} \ar@/_1.7pc/[rr]_-{a} &  \Alb(X)  \ar[r]^-{a'} & h_0^{et}(C)}$$
Indeed, the morphism $a$ induces a morphism from $X$ to $G$. The universal property of the Serre-Albanese scheme gives the morphism $a':\Alb(X)\to G$. The fact that $a=a'\circ \theta$ follows immediately from $\Hom_{\PST(k)}(\Z_{tr}(X),F)=F(X)$ valid for any presheaf with transfers $F$.

By Lemma \ref{lemm-surj-to-alb}, 
the morphism $\theta$ is surjective for the {\'e}tale topology (up to $p$-torsion). This implies that $\cE=\text{Im}(a')$. 
The $1$-generation of $\cE$ follows now from the $1$-generation of $\Alb(X)$.

\smallskip

\noindent
\emph{Step 2:}
By definition we have $\cF\cong h^{\et}_0(\sigma_1^*\cF_1)$ where $\cF_1 = (\sigma_1)_*\cF$ . By Lemma \ref{limn}, $\cF_1$ is a colimit of representable functors:
$$\cF_1=\colim{(C\rightarrow \cF_1)_{\leq 1}} \Z_{\leq 1}(C)$$
with $C$ smooth of dimension $\leq 1$. It follows that 
$$\cF=h_0^{\et}(\sigma_1^*\cF_1)=\colim{(C\rightarrow \cF_1)_{\leq 1}} h_0^{\et}(C)$$
Let $\mathcal{E}\subset \cF$ be a subsheaf. Let's show that $\mathcal{E}$ is $1$-generated. The obvious morphism:
$$\colim{(C \to \cF_1)_{\leq 1}}  \!\! \xymatrix{ h_0^{\et}(C)\times_{\cF} \cE \ar@{->}[r] &  \mathcal{E}}$$
is surjective, even as a presheaf morphism. Indeed, if $\alpha$ is a section of $\mathcal{E}$ over some smooth $k$-variety, there exist objects $(C_i\to \cF_1)_{i=1,\dots, n}$ such that $\alpha$ is in the image of 
$\coprod_{i=1}^n h_0^{\et}(C_i) \to \cF$. Let $C=\coprod_{i=1}^n C_i$. 
Then, $\alpha$ is also in the image of $h_0^{\et}(C)\times_{\cF} \mathcal{E} \to \mathcal{E}$.

Each subsheaf $h_0^{\et}(C)\times_{\cF} \cE \subset h_0^{\et}(C)$ is $1$-generated, by Step 1. This proves that $\cE$ is $1$-generated.

\smallskip

\noindent
\emph{Step 3:}
To finish the proof, we show that any $1$-generated homotopy invariant sheaf is $1$-motivic (proving the first part of Corollary \ref{1genh=1mot} below). Let  $\cF$ be such a sheaf, the surjection $\xymatrix{\sigma_1^*\sigma_{1*}(\cF) \ar@{->>}[r] & \cF}$ factors through $\xymatrix{h_0^{\et}\sigma_1^*\sigma_{1*}(\cF) \ar@{->>}[r] & \cF}$. 
Let $\cN$ be the kernel of the latter surjection.
By Proposition 
\ref{image} and Remark \ref{exampleh}, we know that $h_0^{\et}\sigma_{1}^*\sigma_{1*}(\cF)$ is $1$-motivic. 
By Step 2, $\cN$ is $1$-generated being a subsheaf of the $1$-motivic sheaf $h_0^{\et}\sigma_1^*\sigma_{1*}(\cF)$. By Corollary \ref{tricky-nmot}, this implies that $\cN=0$.
\end{proof}

\begin{cor}
\label{1genh=1mot}
A $1$-generated homotopy invariant {\'e}tale sheaf is $1$-motivic.
Moreover, $\HI_{\leq 1}^{\et}(k)$ is a Serre subcategory of $Shv_{tr}^{\et}(k)$, \ie stable by subobjects, quotients and extensions. In particular, the inclusion  $\HI_{\leq 1}^{\et}(k)\subset Shv_{tr}^{\et}(k)$ is exact. 
\end{cor}

\begin{proof}
The first part was proven in Step 3 of the proof of Lemma \ref{sub1mot}.
The other claims follow easily from Lemma \ref{sub1mot} and Lemma \ref{nmot-coker-ext}. 
\end{proof}

\begin{lemma}\label{newlemmata}
Let $G$ be a semi-abelian group scheme. Let $\cF\subset G$ be an \'etale subsheaf
 with transfers of $G$ such that $\pi_0(\cF)=0$. Then $\cF$ is represented by a closed subgroup of $G$. 
\end{lemma}

\begin{proof}
By Lemma \ref{sub1mot}, we know that $\cF$ is $1$-motivic.
It follows that $\cF$ is a filtered union of images of 
 $\xymatrix{h_0^{\et}(C) \ar[r] & \cF}$ with $C$ a smooth scheme of dimension $\leq 1$ (\cf Step 2 in the proof of Lemma \ref{sub1mot}):
$$\cF = \bigcup_{C\to \cF}\text{Im}(h_0^{\et}(C) \to \cF)$$
As this union is filtered, we have 
$$\cF=\cF^0 = \bigcup_{C\to \cF}\text{Im}(h_0^{\et}(C) \to \cF)^0$$
where $(\dagger)^0$ denotes the kernel of the surjection $\dagger\to \pi_0(\dagger)$. 
One checks immediately that $\text{Im}(h_0^{\et}(C) \to \cF)^0=
\text{Im}(h_0^{\et}(C)^0\to \cF^0)$. 
Now recall \cite[Sect. 3.4]{V} that $h_0^{\et}(C)$ is represented by $\Alb (C)$ so that $h_0^{\et}(C)^0$ is a semi-abelian variety. 
We thus have
$$\cF = \bigcup_{G'\to \cF}\text{Im}(G' \to \cF)^0$$
where the union is taken over maps $G'\to \cF$ with $G'$ a semi-abelian variety. 
Since the image of $G'\to \cF$ is also the image of $G'\to G$ it is then a semi-abelian variety.  
This proves that $\cF$ is the union of the connected subgroups of $G$ contained in $\cF$. As $G$ is Noetherian, any chain of connected subgroups of $G$ is stationary. This proves our claim.
\end{proof}

\begin{defn}
We say that a $1$-motivic sheaf $\cE$ is \emph{finitely generated} if there exist a semi-abelian group scheme $G$ (\ie such that the connected component of the identity $G^0$  is semi-abelian and $\pi_0(G)$ is finitely generated) and a surjection $q:\xymatrix{G\ar@{->>}[r] &  \cE}$.

If moreover $q$ can be chosen so that $\ker(q)$ is finitely generated (as a $1$-motivic sheaf), we say that $\cE$ is \emph{finitely presented} (or constructible).
\end{defn}

\begin{propose}\label{colim-fp-1-mot}
\text{a)} Let $\cE$ be a finitely presented $1$-motivic sheaf. There is a unique and functorial exact sequence
$$0\to L \to G\to \cE\to 0$$
where $G$ is a semi-abelian group scheme and $L$ a lattice (\ie a torsion free and finitely generated $0$-motivic sheaf). 

\text{b)} Let $\cF$ be a $1$-motivic sheaf. Then $\cF$ is a filtered colimit of finitely presented $1$-motivic sheaves. 
\end{propose}

\begin{proof}
For a) choose a presentation 
$$G_1\to G_0 \to \cE\to 0$$
with $G_0$ and $G_1$ semi-abelian group schemes. Denote by $G_1^0$ the connected component of $G_1$ and let $G'=\coker(G_1^0\to G_0)$. Then $G'$ is a semi-abelian group scheme. Moreover, we have a presentation
$$L'\to G' \to \cE\to 0$$
where $L'=G_1/G_1^0=\pi_0(G_1)$.  
Now let $L''$ be the image of $L'$ in $G'$ and $L''_{tor}\subset L''$ its torsion subsheaf.
We define $L=L''/L''_{tor}$ and $G=G'/L''_{tor}$. Then $L$ is a torsion free finitely generated $0$-motivic sheaf, $G$ is a semi-abelian group scheme and
$$0\to L \to G\to \cE\to 0$$
is an exact sequence. The uniqueness and functoriality of this sequence is easy and left to the reader (see also \cite[Prop. 3.2.3]{BK}). 

We now show part b) of the proposition. We divide the proof in two parts.
\smallskip

\noindent
\emph{Part 1:} 
We first consider the case where $\pi_0(\cF)=0$. 
Let $\mathsf{P}(\cF)$ be the category of all morphisms $a:\cE\to \cF$ such that 
\begin{itemize}

\item $\cE$ is a finitely presented $1$-motivic sheaf with $\pi_0(\cE)=0$,

\item $\ker(a)$ is a $0$-motivic sheaf.

\end{itemize}
We will prove that $\mathsf{P}(\cF)$ is filtered and  $$\colim{\cE\to \cF\in \mathsf{P}(\cF)}\cE\simeq \cF$$

For simplicity, we write $\cE/\cF$ an object $\cE\to \cF$ of $\mathsf{P}(\cF)$. If $\cE/\cF$ and $\cE'/\cF$ are two objects in $\mathsf{P}(\cF)$ there is at most one arrow $(\cE/\cF)\to (\cE'/\cF)$. 
Indeed, let $a_1, \, a_2:\cE\to\cE'$.
By the first part of the proposition we can find a commutative diagram
$$\xymatrix{0\ar[r] & L\ar[r] \ar[d] & G \ar[r] \ar[d]^-{b_i} & \cE \ar[r] \ar[d]^-{a_i} & 0\\
0\ar[r] & L'\ar[r] & G' \ar[r] & \cE' \ar[r] & 0}$$
with $L$ and $L'$ lattices and $G$ and $G'$ semi-abelian varieties. Let $H$ be the coequalizer of $b_1$ and $b_2$. Then $G'\to \cF$ factors trough $H\to \cF$. As the kernel of $G'\to \cF$ is $0$-motivic, we deduce that the kernel of $G'\to H$ is also a $0$-motivic sheaf. This happens only when $G'\to H$ is an isogeny. But then $b_1-b_2$ factors through the torsion points of $G'$. This forces  $b_1=b_2$ as $G$ is connected. By a diagram chase, we deduce that $a_1=a_2$.

By the proof of Lemma \ref{newlemmata} 
we know that $\cF$ is a filtered union of images 
$\im(G\to \cF)$ with $G$ a semi-abelian variety. Given such $a:G\to \cF$, $\ker(a)^0$ is a connected subgroup of $G$ by Lemma \ref{newlemmata} (recall that $(\dagger)^0=\ker(\dagger\to \pi_0(\dagger))$). If $G'=G/\ker(a)^0$, the kernel of the morphism $G'\to \cF$ is $0$-motivic. In particular $G'\to \cF\in \mathsf{P}(\cF)$.

To prove that $\mathsf{P}(\cF)$ is filtered, we pick two objects $\cE_1/\cF$ and $\cE_2/\cF$. By the discussion above, we can find $\cG/\cF\in \mathsf{P}(\cF)$ such that $\im(\cG\to \cF)$ contains both $\im(\cE_i\to \cF)$. We reduce then easily to the case where $\im(\cE_1\to\cF)\subset \im(\cE_2\to \cF)$. Let $\cE_3=\cE_1\times_{\cF}\cE_2\subset \cE_1\times \cE_2$. By Lemma \ref{sub1mot}, $\cE_3$ is a $1$-motivic sheaf and 
$\cE_3^0$ is finitely presented as one easily deduce from Lemma \ref{newlemmata}.  
By construction, $\cE_3\to \cE_1$ is surjective and its kernel $\cN$ is contained in $\ker(\cE_1\to \cF)\times \ker(\cE_2\to \cF)$. In particular, it is $0$-motivic. Let $\cE_4=\coker(\cN\to \cE_2)$. Then 
$\cE_4/\cF\in \mathsf{P}(\cF)$ and we have maps $\cE_1\to \cE_4$ and $\cE_2\to \cE_4$. This proves that $\mathsf{P}(\cF)$ is equivalent to a filtred ordered set.

The surjectivity of 
$\colim{\cE\to \cF\in \mathsf{P}(\cF)}\cE\to \cF$ is clear. For injectivity, we use that $\ker(\cE\to \cF)$ is the filtered union of its finitely generated subsheaves $L$ so that $\im(\cE\to \cF)$ is the filtered colimit of the $\cE/L$.

\smallskip

\noindent
\emph{Part 2:}
Now we treat the general case. For $\cL\subset \pi_0(\cF)$ a subsheaf, let $\mathsf{Q}(\cL)$ be the set of finitely generated subsheaves of $\cL$. 
We consider the class $\mathscr{C}$ of functors 
$\cE:\mathsf{Q}(\cL)\to \HI^{\et}_{\leq 1}(k)/\cF$
which assign to $L\in \mathsf{Q}(\cL)$ a morphism 
$\cE(L)\to \cF$ such that:
\begin{enumerate}

\item $\cE(L)^0/\cF^0\in \mathsf{P}(\cF^0)$,

\item $\pi_0(\cE(L))\to \pi_0(\cF)$ is injective and its image is $L$.

\end{enumerate}
We have an obvious notion of isomorphism between functors in 
$\mathscr{C}$ and the isomorphism classes in $\mathscr{C}$ form a set. Given $\cE$ and $\cE'$ defined on 
$\mathsf{Q}(\cL)$ and $\mathsf{Q}(\cL')$, we write $\cE\leq \cE'$ if 
$\cL\subset \cL'$ and the restriction of 
$\cE'$ to $\mathsf{Q}(\cL)$ is isomorphic to $\cE$. 

By Zorn Lemma, we may pick a maximal functor $\cE:\mathsf{Q}(\cL)\to \HI_{\leq 1}^{\et}(k)/\cF$ in $\mathscr{C}$. 
Let us prove that $\cL=\pi_0(\cF)$. Suppose the contrary and let 
$M\subset \pi_0(\cF)$ not contained in $\cL$. 
We may assume that $M/M_0$ is simple (\ie has no proper non-zero subsheaves) where $M_0=M\cap \cL$.

The inverse image of $M$ along $\cF\to \pi_0(\cF)$ is an extension of $M$ by $\cF^0$. This gives an element in $\Ext^1(M,\cF^0)$ (where the $\Ext^1$ is taken in the category of \'etale sheaves). By Step 1, we have 
$$\Ext^1(M,\cF^0)=\colim{\cE^0/\cF^0\in \mathsf{P}(\cF^0)} \Ext^1(M,\cE^0)$$
It follows that we can find $\cE'(M)\to \cF$ such that
\begin{itemize}

\item $\cE'(M)^0/\cF^0\in \mathsf{P}(\cF^0)$,

\item $\pi_0(\cE'(M))\to \pi_0(\cF)$ is injective and its image is $M$,

\item there is a morphism $\cE(M_0)\to \cE'(M)$ over $\cF$.

\end{itemize}
Let $\cL'=\cL+M$.
We define a functor $\cE'$ on $\mathsf{Q}(\cL')$ in the following way.
If $L\subset \cL$ we take $\cE'(L)=\cE(L)$. 
Suppose that $M\subset L$ and let $L_0=L\cap \cL$. We define $\cE_1'(L)$ by the pushout square
$$\xymatrix{\cE(M_0)\ar[r] \ar[d] & \cE'(M)\ar[d] \\
\cE(L_0) \ar[r] & \cE_1'(L)}$$
We then take $\cE'(L)=\cE_1'(L)/\ker(\cE_1'(L)\to \cF)^0$.
For general $L\not \subset \cL$, we let $\cE'(L)$ be the inverse image of $L$ by 
$\cE'(L+M)\to L+M$. 
One easily checks that we have extended the functor $\cE$ to $\mathsf{Q}(\cL')$. This is a contradiction. 

Fix a functor $\cE:\mathsf{Q}(\pi_0(\cF))\to \HI_{\leq 1}^{\et}(k)/\cF$ in $\mathscr{C}$.  
Let $\mathsf{R}(\cF)\subset \mathsf{Q}(\pi_0(\cF))\times \mathsf{P}(\cF^0)$ be the full subcategory whose objects are 
$(L,\cE^0/\cF^0)$ such that 
$\Hom_{\mathsf{P}(\cF^0)}(\cE(L)^0,\cE^0)\neq \emptyset$. 
Given $(L,\cE^0)$ in $\mathsf{R}(\cF)$ we define 
$T(L,\cE^0)=\cE(L)\coprod_{\cE(L)^0}\cE^0$. We get in this way a functor $T:\mathsf{R}(\cF)\to \HI_{\leq 1}^{\et}(k)/\cF$ such that 
$\colim{\mathsf{R}(\cF)} T\simeq \cF$. That $\mathsf{R}(\cF)$ is filtered is clear. The proposition is proven. 
\end{proof}

\begin{cor}
\label{cor-V-contract-1-0}
Let $\cF$ be a $1$-motivic \'etale sheaf. Then, the Voevodsky contraction $\cF_{-1}=\ihom(\G_m,\cF)$ is a torsion free $0$-motivic \'etale sheaf.
\end{cor}

\begin{proof}
It is clear that $\cF_{-1}$ is an \'etale sheaf. Let us show that it is $0$-generated as a presheaf. This is sufficient by Proposition \ref{ogen} and Remark \ref{n-gen-and-topology}.

Suppose we can write $\cF$ as a filtered colimit of $1$-motivic \'etale sheaves $\cF=\colim{\alpha} \cF_{\alpha}$. As $\ihom(\G_m,-)$ commutes with filtered colimits, we need only to show that each $\ihom(\G_m,\cF_{\alpha})$ is $0$-generated. 
By Proposition \ref{colim-fp-1-mot} we may assume that $\cF$ is finitely presented and hence have a short exact sequence:
\begin{equation}
\label{eq-V-contract-1-0}
0 \to L \to G \to \cF \to 0
\end{equation}
with $L$ a lattice and $G$ a semi-abelian group scheme. 

For a smooth $k$-scheme $X$, we have 
a long exact sequence
\begin{equation}
\label{eqn-cf-minus-1}
\xymatrix{0\ar[r] & L(X) \ar[r] & G(X) \ar[r] & \cF(X) \ar[r] & {\rm H}^1_{\et}(X,L) \ar[r] & \dots}
\end{equation}
Let $k_s/k$ be a separable closure with Galois group $G$ and write $X_{k_s}=X\otimes_k k_s$.  
By the Hochschild-Serre spectral sequence 
we have an exact sequence
$$\xymatrix{0\ar[r] & {\rm H}^1(G,{\rm H}^0_{\et}(X_{k_s},L)) \ar[r] & {\rm H}^1_{\et}(X,L) \ar[r] & {\rm H}^0(G, {\rm H}^1_{\et}(X_{k_s},L)) }$$
By \cite[IX, Prop. 3.6 (ii)]{SGA4} we know that ${\rm H}^1_{\et}(X_{k_s},L)=0$ as the restriction of $L$ to $X_{k_s}$ is isomorphic to a direct sum of copies of the constant sheaf $\Z$.
Moreover, ${\rm H}^0_{\et}(X_{k_s},L)={\rm H}^0_{\et}(\underline{\pi}_0(X)\otimes_k k_s)$.
Let $\cK_0$ denote the presheaf on $(Sm/k)_{\leq 0}$ which associates to the spectrum $\Spec(l)$ of an \'etale $k$-algebra $l$ the group
${\rm H}^1(G,{\rm H}^0_{\et}(\Spec(l\otimes_k k_s),L))$. 
If $\cK=\sigma_0^*\cK_0=\underline{\pi}_{0*}\cK_0$, we get from \eqref{eqn-cf-minus-1} an exact sequence of presheaves with transfers:
\begin{equation}
\label{eqn-cf-minus-1-2}
\xymatrix{0 \ar[r] & L \ar[r] & G \ar[r] & \cF \ar[r] & \cK}
\end{equation}
Moreover, as $\cK$ is homotopy invariant, this is an exact sequence of homotopy invariant presheaves with transfers.

The functor $\ihom(\G_m,-)$ is obviously exact on $\PST(k)$. 
Moreover, $\ihom(\G_m,E)=0$ for $E$ a strongly $0$-generated presheaf.
Thus, we obtain from \eqref{eqn-cf-minus-1-2} an isomorphism 
$\ihom(\G_m,G) \simeq  \ihom(\G_m,\cF)$.
It is well know that $\ihom(\G_m,G)$ is a lattice if $G$ is a semi-abelian group scheme. Using that filtered colimit of lattices is torsion free, we get also that $\cF_{-1}$ is torsion free. 
It is also possible to show directly that multiplication by $n$ is injective on $\cF_{-1}$ by noting that it is surjective on $\G_m$ (for the \'etale topology, up to $p$-torsion).  
\end{proof}

\begin{thm}\label{struct-1-mot}
Let $\cF$ be a $1$-motivic sheaf. There exists an exact sequence in 
$\HI_{\leq 1}^{\et}(k)$:
\begin{equation}
\label{eq-struct-1-mot}
\xymatrix{0 \ar[r] & \cN \ar[r] & \cG \ar[r] & \cF\ar[r] & \cK \ar[r] & 0}
\end{equation}
such that:
\begin{enumerate}
\item[(i)] $\cN$ and $\cK$ are $0$-motivic sheaves, $\cK=\pi_0(\cF)$ and $\cN$ is torsion free.

\item[(ii)] We have an isomorphism
$\ihom(\G_m,\cG)\simeq \ihom(\G_m,\cF)$
\end{enumerate}
Let $\cL=\ihom(\G_m,\cF)$. Then $\cL$ is a torsion free $0$-motivic sheaf and the canonical morphism $\cL\otimes \G_m \to \cG$ is injective. 
Let $\cA=\coker(\cL\otimes \G_m\to \cG)$.

\begin{enumerate}
\item[(iii)] $\cA$ is a filtered colimit of abelian varieties.

\item[(iv)] With rational coefficients, $\cA$ is isomorphic to a direct sum of simple abelian varieties, \ie $\cA\otimes \Q \simeq \oplus_{\beta} B_{\beta}\otimes \Q$. 

\end{enumerate} 

\end{thm}

\begin{proof}
We know by Proposition \ref{colim-fp-1-mot} that $\cF$ is a filtered colimit of finitely presented $1$-motivic sheaves. We get \eqref{eq-struct-1-mot} by taking the colimit of the functorial exact sequences in Proposition \ref{colim-fp-1-mot}.

Let us check the properties (i)-(iv). 
We get (i) by construction. To check that $\ihom(\G_m,\cG)\simeq \ihom(\G_m,\cF)$ we may assume that $\cF$ is finitely presented. Then the claim follows from the proof of  
Corollary \ref{cor-V-contract-1-0}. 

Also, to show that $\cL\otimes \G_m\to \cG$ is injective, we may assume that $\cF$ is finitely presented. Here again, the claim follows from the proof of 
Corollary \ref{cor-V-contract-1-0}. Property (iii) is clear from
Proposition \ref{colim-fp-1-mot}.

It remains to prove (iv). Let $\cB\subset \cA$ be a maximal subsheaf of $\cA$
that can be written as a direct sum of simple abelian varieties (after tensoring by $\Q$). This exists by Zorn Lemma. Assume that $\cB\neq \cA$. By (iii) there is and abelian variety $C$ and $C\to \cA$ whose image is not contained in $\cB$. Dividing by the connected component of the kernel of $C\to \cA$, we may assume that $C \to \cA$ is injective (as the kernel is torsion). Consider now $C\cap \cB$. This is a connected subgroup of $C$. Let $C'$ be a supplement of $C\cap \cB$ in $C$. Then $\cB\oplus C'\subset \cA$. This is a contradiction. 
\end{proof}

\begin{propose} 
\label{1mot-left-adj}
The embedding $\HI_{\leq 1} ^{\et}(k)\into Shv_{tr}^{\et}(k)$ has a  left adjoint  $\Alb : Shv_{tr}^{\et}(k)\to \HI_{\leq 1}^{\et}(k)$ given by the following
$$  \Alb (\cF)\df \colim{X\to \cF} \Alb (X)$$
\end{propose}

\begin{proof}  Let $\cF\in Shv_{tr}^{\et}(k)$ and  $\cE\in  \HI_{\leq 1}^{\et}(k)$. Consider the following commutative diagram, \cf Lemma \ref{lim}:
$$\xymatrix{ \Hom (\Alb (\cF), \cE)    \ar[r]\ar@{=}[d] & \Hom(\cF,\cE)\ar@{=}[d] \\
\Hom ( \underset{X \to \cF}{\text{Colim}} \; \Alb (X), \cE)\ar[r]\ar@{=}[d] &\Hom (\underset{X \to \cF}{\text{Colim}} \; \Z_{tr}(X), \cE)\ar@{=}[d] \\
\underset{X \to \cF}{\text{Lim}} \; \Hom (\Alb (X), \cE)\ar[r]&\underset{X \to \cF}{\text{Lim}} \; \Hom (\Z_{tr}(X), \cE)}$$
We are then left to show the following:
\end{proof}

\begin{lemma}\label{adjalb} For $\cE\in  \HI_{\leq 1}^{\et}(k)$
$$\alpha_X:\xymatrix{\Hom (\Alb (X), \cE)
 \ar[r]^-{\sim} & \Hom (\Z_{tr}(X), \cE)}$$
is invertible.
\end{lemma}

\begin{proof}
We may assume $k$ separably closed by Lemma \ref{hom-Galois}
below.
As $\xymatrix{\Z_{tr}(X) \ar[r] & \Alb(X)}$ is a surjection of {\'e}tale sheaves (again up to $p$-torsion), our homomorphism is injective. We only need to check that $\alpha_X$ is surjective. 
Take $s\in \Hom(\Z_{tr}(X),\cE)$. 

By Proposition \ref{colim-fp-1-mot}, we know that $\cE$ is a filtered colimit of finitely presented $1$-motivic sheaves. Since {\'e}tale topology is quasi-compact,  $s$ factors through  $\cE_0\to \cE$ with $\cE_0$ finitely presented. We may then assume $\cE$ to be itself finitely presented. 
We then have an exact sequence
$$\xymatrix{0 \ar[r] & L \ar[r] & G \ar[r] & \cE \ar[r] & 0}$$
with $L$ a lattice and $G$ a semi-abelian group scheme. We deduce a long exact sequence in cohomology:
$$\xymatrix{0 \ar[r] & \Hom(\Z_{tr}(X) , L) \ar[r] & \Hom(\Z_{tr}(X),G) \ar[r] & \Hom(\Z_{tr}(X), \cE) }$$
\hfill $\xymatrix{&\ar[r] & \text{H}_{\et}^1(X,L) \ar[r] & \dots}$

As $k$ is separably closed, $L$ is isomorphic to the constant sheaf $\Z^r$.
By \cite[IX, Prop. 3.6 (ii)]{SGA4}, $\text{H}^1_{\et}(X,\Z^r)=0$ since $X$ is smooth and hence normal.
It follows that $s$ factors:
$$\xymatrix{\Z_{tr}(X) \ar[r] \ar@/_1.5pc/[rr]_-s & G \ar[r] & \cE}$$
By the universality of $\Alb(X)$ we get a further factorization:
$$\xymatrix{\Z_{tr}(X) \ar[r] \ar@/_1.7pc/[rrr]_-s  & \Alb(X) \ar[r] \ar@/^1.5pc/[rr]^-{s'} & G \ar[r] & \cE }$$ 
Then $s'$ is mapped to $s$ by $\alpha_X$. This proves the surjectivity of $\alpha_X$.
\end{proof}

\begin{lemma}
\label{hom-Galois}
Let $\cA$ and $\cB$ be two \'etale sheaves with transfers on $Sm/k$. Let $k\subset k'$ be a Galois extension and denote by $\cA_{k'}$ and $\cB_{k'}$ the pull-backs to $Sm/k'$. Then we have an isomorphism:
$$\xymatrix{\Hom(\cA,\cB) \ar[r]^-{\sim} & \Hom(\cA_{k'}, \cB_{k'})^{\Gal(k'/k)}}$$ 
\end{lemma}

\begin{defn}
\label{1mot}
Denote by $(-)^{\leq 1}:\HI_{tr}^{\tau}(k) \to \HI_{\leq 1}^{\tau}(k)$ the restriction of $\Alb$ to $\HI_{tr}^{\tau}(k) \subset Shv_{tr}^{\tau}(k)$. It is left adjoint to the inclusion $\HI_{\leq 1}^{\tau}(k) \subset \HI_{tr}^{\tau}(k) $.
\end{defn}

\begin{remark}
 The category $\HI_{\leq 1}^{\et}(k)$ is the smallest cocomplete Serre subcategory of $Shv_{tr}^{\et}(k)$ containing lattices and {\'e}tale sheaves represented by semi-abelian varieties. It is also the smallest cocomplete Serre subcategory containing $h_0^{\et}(C)$ for $C$ smooth curves. 
\end{remark}

\subsection{$n$-generated for $n\geq 2$} In this section, we propose a
conjecture  that makes it possible to extend the results about
$0$-motivic and $1$-motivic sheaves to the $n$-motivic case. Here also
assume that $p$ is inverted and take $\tau=\text{{\'e}t}$.

\begin{conj}
\label{weak-BBF}
For any smooth $k$-variety $X$, there exists a filtration $F^{i+1}h_0^{\et}(X) \subset F^ih_0^{\et}(X)$ such that:
\begin{enumerate}

\item[(A)] $F^0h_0^{\et}(X)=h_0^{\et}(X)$ and $F^nh_0^{\et}(X)=0$ for $n\geq \dim(X)+1$,

\item[(B)] The filtration is compatible with the action of correspondences, \ie for $\gamma\in Cor(X,Y)$ the induced morphism of homotopy sheaves $h_0^{\et}(X) \longby{\gamma} h_0^{\et}(Y)$ is compatible with the filtration.

\item[(C)] If $U$ is a dense open subvariety of $X$ then $h_0^{\et}(U) \to h_0^{\et}(X)$ is strict for the filtration. 

\item[(D)] For $n\geq 0$, the quotient $F^0h_0^{\et}(X)/F^{n+1}h_0^{\et}(X)$ is $n$-generated.  
\end{enumerate} 
\end{conj}

\begin{remark}
When $X$ is smooth projective,  the Bloch-Beilinson conjectural filtration on the Chow group of $0$-cycles induces a filtration on $h_0(X)_{\Q}$ as we have $h_0(X)_{\Q}(K)=\text{CH}_0(X\times_k K)_{\Q}$ for any $k$-field $K$. This filtration should be the same as the one predicted in \ref{weak-BBF}. We remark also that the properties of the Bloch-Beilinson filtration imply (A) and (B) in the case that $X$ is projective (at least with rational coefficients). Moreover, with more effort, one should obtain (D) as well.
\end{remark}

\begin{lemma}\label{subnmot}
Assume (A), (B) and (D) of Conjecture \ref{weak-BBF}.
Let $\cF$ be an $n$-motivic sheaf. Any subsheaf of $\cF$ is again $n$-motivic.
\end{lemma} 

\begin{proof}
The proof is very similar to Lemma \ref{sub1mot}. One argues in three steps. The second and third steps are formal and extend literally to the general case. The first step is to show that any subsheaf of $\cF=h_0^{\et}(X)$ is $n$-generated if $\dim(X)\leq n$. As in \ref{sub1mot}, we can suppose that $\cE=\text{Im}(a:h_0^{\et}(W) \to h_0^{\et}(X))$ for some smooth variety $W$. 
As $F^{n+1}h_0^{\et}(X)=0$ and $a$ is compatible with the filtration of \ref{weak-BBF} we get a factorization:
$$\xymatrix{h_0^{\et}(W) \ar@{->>}[r] \ar@/_1.9pc/[rr]_-a & h_0^{\et}(W)/F^{n+1} h_0^{\et}(W) \ar[r]^-{a'} & h_0^{\et}(X)  }$$
It is clear that the image of $a$ is equal to the image of $a'$. This proves that $\cE$ is $n$-generated since it is a quotient of the $n$-generated sheaf $h_0^{\et}(W)/F^{n+1} h_0^{\et}(W)$.
\end{proof}

\subsubsection{} It follows from Lemma \ref{subnmot} that under (A) and (B) of Conjecture \ref{weak-BBF}, condition (D) is equivalent to the stronger one:
\begin{enumerate}

\item[($\text{D}^{\prime}$)]  For $n\geq 0$, the quotient $F^0h_0^{\et}(X)/F^{n+1}h_0^{\et}(X)$ is $n$-motivic.

\end{enumerate}

We get also the following (\cf \ref{1genh=1mot} for $n=1$):

\begin{cor}
\label{nmot-thick}
Assume (A), (B) and (D) of Conjecture \ref{weak-BBF}. Then $\HI_{\leq
  n}^{\et}(k)$ is a Serre abelian subcategory of $\HI_{tr}^{\et}(k)$.
  Moreover, the inclusion $\HI_{\leq n}^{\et}(k)\subset \HI_{tr}^{\et}(k)$ is exact.
\end{cor}

\begin{proof}
We just saw in Lemma \ref{subnmot} that $\HI_{\leq n}^{\et}(k)$ is stable by subobjects. Stability by cokernels and extensions is proven in  Lemma \ref{nmot-coker-ext}.
\end{proof}

Another consequence of Conjecture \ref{weak-BBF} is the following (\cf \ref{1mot-left-adj} for $n=1$):

\begin{propose}\label{n-conjecture}
Assume Conjecture \ref{weak-BBF} with rational coefficients. There exist left adjoints 
$$(-)^{\leq n}:\xymatrix{\HI^{\et}_{tr}(k)_{\Q} \ar[r] & \HI^{\et}_{\leq n}(k)_{\Q}}$$ 
to the inclusions $\HI^{\et}_{\leq n}(k)_{\Q}\subset \HI^{\et}_{tr}(k)_{\Q}$.
Moreover, the functors $(-)^{\leq *}$ and the filtration $F^*$ are related by the following:
$$  (\cF)^{\leq n}\cong \colim{X\to \cF} h^{\et}_0(X)_{\Q}/F^{n+1}h^{\et}_0(X)_{\Q}$$
and $F^nh^{\et}_0(X)_{\Q}=\emph{Ker}(h^{\et}_0(X)_{\Q} \to (h^{\et}_0(X)_{\Q})^{\leq n-1})$.

Conversely, if the adjoints $(-)^{\leq n}$ exist and the $\HI^{\et}_{\leq n}(k)_{\Q}$ are Serre subcategories of $\HI^{\et}_{tr}(k)_{\Q}$  for all $n\geq 0$, then Conjecture \ref{weak-BBF} holds.
\end{propose}

\begin{proof} First assume Conjecture \ref{weak-BBF}.  As in the $1$-motivic case, for $\cE\in  \HI^{\et}_{\leq n}(k)_{\Q}$ we are left to show that
$$\alpha_X:\xymatrix{\Hom (h_0^{\et}(X)^{\leq n}_{\Q}, \cE)\ar[r] &  \Hom (\Q_{tr}(X), \cE)}$$
is surjective (note that injectivity is clear). As in the proof of Lemma \ref{adjalb} we may suppose that $\cE$ is the quotient of a $h^{\et}_0(Y)_{\Q}$ for $Y$ of dimension $\leq n$. Since we are working with rational coefficients the map $\xymatrix{h^{\et}_0(Y)_{\Q}\ar@{->>}[r] & \cE}$ is a surjection of Zariski sheaves. It follows that for $s\in \text{H}^0(X, \cE)= \Hom (\Q_{tr}(X), \cE)$ there exists a dense open subset $U$ of $X$ such that $s_{\mid U}$ lifts to $h_0^{\et}(Y)_{\Q}$:
$$\xymatrix{h_0^{\et}(Y)_{\Q} \ar@{->>}[r] & \cE \\
\Q_{tr}(U) \ar[r] \ar@{.>}[u]^t & \Q_{tr}(X)\ar[u]_-s }$$
As $\cE$ and $h^{\et}_0(Y)_{\Q}$ are homotopy invariant $s$ and $t$ factors through $h^{\et}_0(X)_{\Q}$ and $h^{\et}_0(U)_{\Q}$:
$$\xymatrix{h^{\et}_0(Y)_{\Q} \ar@{->>}[r] & \cE \\
h^{\et}_0(U)_{\Q} \ar[r]_{u} \ar[u]^{t_0} & h^{\et}_0(X)_{\Q}\ar[u]_{s_0} }$$
By Conjecture \ref{weak-BBF}, the map $t_0$ is compatible with the filtration. It sends the subsheaf $F^{n+1}h^{\et}_0(U)_{\Q}$ to $F^{n+1}h_0^{\et}(Y)_{\Q}=0$. 

The morphism $u$ is surjective. To see this, it suffices by Yoneda to show that 
$\Hom(h_0^{\et}(X)_{\Q},\dagger)\to \Hom(h_0^{\et}(U)_{\Q},\dagger)$ is injective for any homotopy invariant \'etale sheaf of $\Q$-vector spaces $\dagger$. This map is nothing but $\dagger(X)\to \dagger(U)$ which is injective by \cite[Lemma 22.8]{VL}.

By Conjecture \ref{weak-BBF} (C), $\xymatrix{F^{n+1}h_0^{\et}(U)_{\Q} \ar@{->>}[r] & F^{n+1}h_0^{\et}(X)_{\Q}}$ is surjective. This implies that $s_0$ maps $F^{n+1}h^{\et}_0(X)_{\Q}$ to $0$. This gives a factorization:
$$\xymatrix{\Q_{tr}(X) \ar[r] \ar@/^1.7pc/[rrr]^-{s} & h^{\et}_0(X)_{\Q} \ar[r] \ar@/_1.5pc/[rr]_-{s_0} & h^{\et}_0(X)_{\Q}/F^{n+1}h^{\et}_0(X)_{\Q} \ar@{.>}[r] & \cE}$$
The dotted arrow is mapped to $s$ by $\alpha_X$.

\smallskip

Conversely, suppose that the left adjoints $(-)^{\leq n}$ exist for all $n\geq 0$ and define $F^n$ as in the statement for any $\cF\in \HI^{\et}_{tr}(k)_{\Q}$ to be the kernel of $\cF \to (\cF)^{\leq n-1}$. Properties (A), (B) and (D) are clear. We need only to check (C). 

First remark that the inclusion $\HI^{\et}_{\leq n}(k)_{\Q} \subset
\HI^{\et}_{tr}(k)_{\Q} $ is exact (as it admits a left adjoint). We will prove more generally that for any surjective morphism:
$$a:\xymatrix{\cE \ar@{->>}[r] & \cF}$$
the induced morphism $F^{n+1}(a):\xymatrix{F^{n+1}(\cE) \to F^{n+1}(\cF)}$ is again surjective. Let us denote by $\cK$ the cokernel of $F^{n+1}(a)$. It is sufficient to prove that $\cK$ is $n$-motivic. Indeed, in this case the cokernel $\cL$ of $F^{n+1}(\cE) \to \cF$ is 
$n$-motivic being an extension 
of two $n$-motivic sheaves:
$$\xymatrix{0 \ar[r] & \cK \ar[r] & \cL \ar[r]^-{\dag} & (\cF)^{\leq n} \ar[r] & 0}$$ 
The universality of $(\cF)^{\leq n}$ implies that $\dag$ is invertible. This forces $\cK$ to be zero.

To check that $\cK$ is $n$-motivic, consider the diagram:
$$\xymatrix{&  &  & \cN \ar[d]  & \\
0 \ar[r] & F^{n+1}\cE \ar[r] \ar[d] & \cE \ar[r] \ar@{->>}[d] & (\cE)^{\leq n} \ar[r] \ar@{->>}[d] & 0 \\
0 \ar[r] & F^{n+1}\cF \ar[r] \ar@{->>}[d]  & \cF \ar[r] \ar[d] & (\cF)^{\leq n} \ar[r]  & 0 \\
& \cK \ar[r] & 0 & }$$
where $\cN$ is the kernel of $(\cE)^{\leq n} \to (\cF)^{\leq n}$ which is $n$-motivic. By the snake lemma $\cK$ is a quotient of $\cN$. If we further assume that $\HI_{\leq n}(k)_{\Q} \subset \HI_{tr}(k)_{\Q}$ is a Serre subcategory, $\cK$ is even $n$-motivic.
\end{proof}

\begin{remark}
Proposition \ref{n-conjecture} shows that if a filtration $F^i$ as in Conjecture \ref{weak-BBF} exists then it is unique (at least after tensoring with $\Q$). 
\end{remark}

\section{Deriving $\pi_0$ and $\Alb$}

\subsection{Generalities} We first explain a general technique to derive right exact functors between Grothendieck abelian categories.
For an abelian category $\catA$, denote by $C(\catA)$ the category of complexes of objects of $\catA$,  $K(\catA)$ the homotopy category of $C(\catA)$ and $D(\catA)$ the
derived category of $\catA$. When $\catA$ is Grothendieck, by a theorem of Joyal (\cf \cite{Joyal, beke}) the category $C(\catA)$  has a model category structure where
the \emph{cofibrations} are the injective morphisms and the \emph{weak
  equivalences} are the quasi-isomorphisms. In particular $D(\catA)$
exists without enlarging the universe (see also the remark of Gabber in \cite[10.4.5]{Weibel}).  In the sequel we will use the homological indexing for complexes. 

\begin{lemma}
Let $\catA$ be a Grothendieck abelian category and $I_{\d}\in
C(\catA)$ a fibrant complex. For any  $A_{\d}\in C(\catA)$ we have an isomorphism: $\Hom_{K(\catA) }(A_{\d},I_{\d})\simeq \Hom_{D(\catA)}(A_{\d},I_{\d})$.
\end{lemma}

\begin{proof}
As $A_{\d}$ is cofibrant and $I_{\d}$ is fibrant we know that
$\Hom_{D(\catA)}(A_{\d},I_{\d})$ coincides with the homotopy classes of
maps in $\Hom_{C(\catA)}(A_{\d},I_{\d})$ with respect to a fixed
cylinder $\text{Cyl}_{A_{\d}}$ (see \cite[II.1]{Jar-Groess}). 
When we take the cylinder to be the
cone of $(id,-id):A_{\d} \to A_{\d}\oplus A_{\d}$ we get the usual
homotopy relation on maps of complexes.
\end{proof}

\begin{remark}
Let $\catA$ be a Grothendieck abelian category and $I_{\d}\in
C(\catA)$ a fibrant complex. For any $n\in \Z$ the object $I_n$ is
injective. Indeed, we may assume $n=0$. Let $A\subset B$ and fix $A\to
I_0$. We denote by $N$ the kernel of the composition $A\to I_0\to
I_1$. We get then a morphism of complexes:
$$\xymatrix{ & 0 \ar[r] \ar[d] & A \ar[r] \ar[d] & A/N \ar[r] \ar[d] & 0
  \ar[d] &  \\
\dots \ar[r] & I_{-1} \ar[r] & I_0 \ar[r] & I_1 \ar[r] & I_2 \ar[r] &
  \dots }$$
Using the left lifting property of $I_{\d} \to 0$ with respect to the
  trivial cofibration: 
$$\xymatrix{  0 \ar[r]  & A \ar[r] \ar[d] & A/N \ar[r] \ar[d] & 0
    \\
0 \ar[r]  & B \ar[r]  & B/N \ar[r] & 0 }$$
we get an extension $B \to I_0$ of $A\to I_0$. 
\end{remark}

\begin{remark}
\label{rem-G(fib)=fib}
Let $F:\catA \to \catB$ be an exact functor between Grothendieck abelian categories and suppose that $G:\catB \to \catA$ is right adjoint to $F$.
Then:
$$(F,G):\xymatrix{C(\catA) \ar[r] & C(\catB)}$$
is a Quillen adjunction for the Joyal model structures. Indeed, $F$ preserves cofibrations and quasi-isomorphisms. In particular,
$G$ takes fibrant complexes to fibrant complexes. 
\end{remark}

\subsubsection{} Any left exact functor $G:\catB \to \catA$ between Grothendieck abelian categories admits a total right derived functor:
$$\R G:D(\catB) \to D(\catA)$$
Suppose that $G$ admits a left adjoint $F:\catA \to
  \catB$. We will describe a condition (see Proposition \ref{existance-LF}) which grants the existence of
a total left derived functor $\L F$ left adjoint to $\R G$. 
This condition is directly inspired from \cite[Def. 1.49]{MV}. From now
on, we implicitly assume our abelian categories to be Grothendieck.

\begin{defn}
A complex $P_{\d}\in C(\catA)$ is $F$-\emph{admissible} if for any fibrant complex $I_{\d}\in C(\catB)$ we have an isomorphism: $$\Hom_{K(\catA)}(P_{\d},G(I_{\d}))\simeq \Hom_{D(\catA)}(P_{\d},G(I_{\d}))$$
\end{defn}

\begin{propose} 
\label{existance-LF}
If $C(\catA)$ has enough $F$-admissible complexes, \ie any $A_{\d}\in C (\catA)$ is quasi-isomorphic to an $F$-admissible complex, then $F$ admits a total left derived functor $\L F:\xymatrix{D(\catA) \ar[r] & D(\catB)}$. Furthermore, $\L F$ is a left adjoint of $\R G$. 
\end{propose}

\begin{proof}
Let $A_{\d}\in C(\catA)$ and $B_{\d}\in C(\catB)$. Choose quasi-isomorphisms $P_{\d}\simeq A_{\d}$ and $B_{\d}\simeq I_{\d}$ with $P_{\d}$ $F$-admissible and $I_{\d}$ fibrant.
We then get isomorphisms:
$\Hom_{D(\catA)}(A_{\d}, \text{R} G (B_{\d})) \simeq \Hom_{D(\catA)}(P_{\d}, G(I_{\d}))\simeq \Hom_{K(\catA)}(P_{\d},G(I_{\d})) 
 \simeq \Hom_{K(\catB)}(F(P_{\d}),I_{\d})\simeq \Hom_{D(\catB)}(F(P_{\d}),B_{\d})$. This shows that the covariant functor $\Hom_{D(\catA)}(A_{\d},\R G(-))$ is co-represented by $F(P_{\d})$. This proves the existence of a left adjoint to $\R G$.
\end{proof}

\subsubsection{} We give some lemmas that help in proving the
existence of enough $F$-admissible complexes for a Grothendieck abelian category. The following is a direct analogue of the second statement in \cite[Lemma 1.53]{MV}.

\begin{lemma}
\label{F-adm-colim}
The full subcategory of $C(\catA)$ whose objects are the $F$-admissible complexes is stable by cones and arbitrary sums.
Furthermore, suppose given a diagram:
$$ \xymatrix{(P_{o})_{\d} \ar[r]^-{a_0} & (P_1)_{\d} \ar[r]^-{a_1} & \dots \ar[r] & (P_n)_{\d} \ar[r]^-{a_n} & \dots   }$$
 of $F$-admissible complexes such that $a_n$ and $F(a_n)$ are
 injective for all $n\geq 0$. Then the colimit $P_{\d}$ (computed in $C(\catA)$) of the above diagram is again $F$-admissible. 
\end{lemma}

\begin{proof}
Only the last statement needs a proof. 
Let $I_{\d}\in C(\catB)$ be fibrant and choose a fibrant replacement
$G(I_{\d})\simeq J_{\d}$. As usual we denote $\Hom_{\d}$ the total complex associated to the double complex of degreewise morphisms of chain complexes.
We then have two isomorphisms:
\begin{equation}
\label{F-adm-colim-1}
\xymatrix{\Hom_{\d}(P_{\d},G(I_{\d}))  \ar[r]^-{\sim} & \lim{n} \Hom_{\d}((P_n)_{\d},G(I_{\d})) }
\end{equation}
\begin{equation}
\label{F-adm-colim-2}
\xymatrix{\Hom_{\d}(P_{\d},J_{\d})  \ar[r]^-{\sim} & \lim{n} \Hom_{\d}((P_n)_{\d},J_{\d}) }
\end{equation}
We know by hypothesis that
$\Hom_{\d}((P_n)_{\d},G(I_{\d})) \to \Hom_{\d}((P_n)_{\d},J_{\d})$ is
a quasi-isomorphism for all $n$. In order to conclude, we need to know
that the limits in \eqref{F-adm-colim-1} and \eqref{F-adm-colim-2} are
actually homotopy colimits. This follows from the fact that
$$\Hom_{\d}(F((P_n)_{\d}),I_{\d}) \to
\Hom_{\d}(F((P_{n-1})_{\d}),I_{\d})$$
and 
$$\Hom_{\d}((P_n)_{\d},J_{\d}) \to
\Hom_{\d}((P_{n-1})_{\d},J_{\d})$$
are surjective as $I_{\d}$ and $J_{\d}$ are componentwise injective.
\end{proof}

\begin{cor}
\label{adm-generator}
Keep the notation as above. Suppose we have a generator $E\in \catA$ which is $F$-admissible as a complex concentrated in degree $0$. Then $C(\catA)$ has enough $F$-admissible complexes.  
\end{cor}

\begin{proof}
For an object $A\in \catA$, we define a complex $P(A)_{\d}$ with a quasi-isomorphism $P(A)_{\d} \to A$ such that: 
\begin{itemize}

\item  $P(A)_{n}=0$ for $n<0$,  

\item  $P(A)_0=\coprod_{E\to A} E$,

\item and for $n>0$:
$$P(A)_n=\coprod_{E \to {\rm ker}(P(A)_{n-1} \to P(A)_{n-2})} E$$

\end{itemize}
The complex $P(A)_{\d}$ is functorial in $A$. We define $P(K)_{\d}$ for a bounded complex $K=K_{\d}$ as the simple complex associated to $P(K_{\d})_{\d}$. If $K_{\d} \to L_{\d}$ is a monomorphism of complexes, $P(K)_{\d}\to P(L)_{\d}$ is then a split monomorphism in each degree. In particular, $F(P(K)_{\d})\to F(P(L)_{\d})$ is injective.

Now, let $K=K_{\d}\in C(\catA)$. We may write 
$K=\colim{n\in \N} \tau_{\geq -n} \sigma_{\leq n}(K)$
where, $\tau_{\leq -n}$ is the good truncation and $\sigma_{\leq n}$ is the bad one. We will show that $\colim{n\in \N} P(\tau_{\geq -n} \sigma_{\leq n}(K))_{\d}$
is $F$-admissible. By the last statement of
\ref{F-adm-colim}, it suffices to show that each $P(\tau_{\geq -n} \sigma_{\leq n}(K))_{\d}$ is $F$-admissible. We are thus reduced to the case where $K_{\d}$ is bounded. Using the stability of $F$-admissibility by mapping cone (\cf \ref{F-adm-colim}), we may further suppose that $K=A[0]$ is concentrated in degree zero.

To show that the complex $P(A)_{\d}$ is $F$-admissible, we write it as the colimit of $\sigma_{\leq n} (P(A)_{\d})$ and use again \ref{F-adm-colim}.
\end{proof}

\begin{lemma}
\label{lemma-reduction-inj}
Let $P\in \catA$ such that
$\Hom_{D(\catA)}(P[0],-)$ commutes with arbitrary sums. Then $P$ is
$F$-admissible if and only if for any injective $I\in \catB$ we have
$\Ext^{i}(P,G(I))=0$ for $i>0$.
\end{lemma}

\begin{proof}
The condition is clearly necessarily as $I[0]$ is a fibrant
complex. Let us show that it is sufficient. 

For a fibrant complex $I_{\d}$ in $C(\catB)$ one has
$$I_{\d}=\colim{n} (\holim{m} \sigma_{\leq -n}( \sigma_{\geq m}(I_{\d})))$$
with $\sigma_{\leq -n}$ and $\sigma_{\geq m}$ the bad truncations of complexes.
We are then reduced to the case where $I_{\d}$ is a bounded complex of
injective objects. By induction we might further assume $I_{\d}$ concentrated in one degree. That $\Hom_{K(\mathbf{A})}(P[0],G(I)[n])\to \Hom_{D(\mathbf{A})}(P[0],G(I)[n])$ is invertible is clear if $n<0$ and follows from our assumption of $n>0$.
\end{proof}

\subsection{The functors $\L\sigma_n^*$} As an application we get:

\begin{lemma}\label{derupstar} 
The functor  
$$\sigma_{n*}=\R\sigma_{n*}: D(Shv_{tr}^{\tau}(k))\longby{} D(Shv^{\tau}_{tr}(k_{\leq n}))$$ 
has a left adjoint $\L\sigma_n^*: D(Shv_{tr }^{\tau}(k_{\leq n}))\longby{ }D(Shv_{tr}^{\tau}(k))$.
\end{lemma}

\begin{proof} 
We need to check the existence of enough $\sigma_n^*$-admissible complexes in $C(Shv_{tr}^{\tau}(k_{\leq n}))$.
By Lemma \ref{adm-generator} it is sufficient to prove that for any smooth $k$-variety $X$ of dimension $\leq n$, the complex concentrated in degree zero $\Z_{\leq n}(X)$ is $\sigma_n^*$-admissible.

Let $\cI_{\d}$ be a fibrant complex in $C(Shv_{tr}^{\tau}(k))$ and choose a fibrant resolution $\xymatrix{\sigma_{n*} \cI_{\d} \ar[r] & \cJ_{\d}}$. By the commutative diagram:
 $$\xymatrix{\Hom_{K(Shv_{tr}^{\tau}(k_{\leq n}))}(\Z_{\leq n}(X), \sigma_{n*}\cI_{\d}) \ar[r] \ar@{=}[d] &  \Hom_{K(Shv_{tr}^{\tau}(k_{\leq n}))}(\Z_{\leq n}(X),\cJ_{\d})\\
\Hom_{K(Shv_{tr}^{\tau}(k))}(\Z_{tr}(X), \cI_{\d}) \ar@/_/[ur]_-{a} &  }$$ 
We need to show that $a$ is invertible. But by Lemma \ref{teck-exist-sigma-n-adm} below we have:
$$\Hom_{K(Shv_{tr}^{\tau}(k))}(\Z_{tr}(X), \cI_{\d})\simeq \HH^0(X,(\cI_{\d})_{\mid X_{\tau}})$$
and also (see Remark \ref{rem-idiote})
$$\Hom_{K(Shv_{tr}^{\tau}(k_{\leq n}))}(\Z_{\leq n}(X),\cJ_{\d}) \simeq \HH^0(X,(\cJ_{\d})_{\mid X_{\tau}})$$
where $X_{\tau}$ is the category $\text{\'Et}/X$ of $X$-{\'e}tale schemes together
with the $\tau$-topology.

The result follows then from the fact that $(\cI_{\d})_{\mid X_{\tau}}$
is quasi-isomorphic to $(\cJ_{\d})_{\mid X_{\tau}}$.
\end{proof}

\begin{lemma}\label{teck-exist-sigma-n-adm}
Let $\cI_{\d}\in C(Shv_{tr}^{\tau}(k))$ be a fibrant complex. Then 
$$\Hom_{K(Shv^{\tau}_{tr}(k))}(\Z_{tr}(X), \cI_{\d}) \simeq \HH^0(X, (\cI_{\d})_{\mid X_{\tau}})$$
with $\mathbb{H}^*(X,-)$ the $\tau$-hypercohomology of $X$.
\end{lemma}

\begin{proof}
This is due to Voevodsky. Let us recall quickly his proof. 
We may assume $\tau \in \{{\rm Nis},\text{\'et}\}$. The Nisnevich and \'etale cohomology can be computed using \v{C}ech hypercovers. Giving a $\tau$-cover $f:\xymatrix{X' \ar[r] & X}$ by an \'etale morphism, we need to show that:
\begin{equation}
\label{eq-lemme-tau-fib}
\xymatrix{\Gamma(X,\cI_{\d}) \ar[r] & \Gamma(\check{C}(f),\cI_{\d})}
\end{equation}
is a quasi-isomorphism (where $\check{C}(f)$ is the \v{C}ech hypercover associated to $f$). The morphism \eqref{eq-lemme-tau-fib} is equal by adjunction to:
$$\xymatrix{\Hom_{\d}(\Z_{tr}(X),\cI) \ar[r] & \Hom_{\d}(\Z_{tr}(\check{C}(f)),\cI)}$$
As $\cI$ is fibrant, we only need to show that $\Z_{tr}(\check{C}(f)) \to \Z_{tr}(X)$ is a quasi-isomorphism of complexes of $\tau$-sheaves. This is true by \cite[Prop. 6.12]{VL}.
\end{proof}

\begin{remark}
\label{rem-idiote}
The statement of Lemma \ref{teck-exist-sigma-n-adm} holds for $Shv^{\tau}_{tr}(k_{\leq n})$. The same proof works with obvious changes.
\end{remark}

\begin{lemma}
\label{Lsigma-up-star}
The unit of the adjunction $\xymatrix{{\rm id} \ar[r]^-{\sim} & \R \sigma_{n*} \L\sigma_n^*}$ is invertible.
\end{lemma}

\begin{proof}
As a triangulated category with arbitrary sums $D(Shv_{tr}^{\tau}(k_{\leq n}))$ is generated by $\Z_{\leq n}(X)[0]$ for $X\in (Sm/k)_{\leq n}$. 
As both $\R\sigma_{n*}$ and $\L\sigma_n^*$ commute with arbitrary sums, we only need to prove that:
$$\xymatrix{\Z_{\leq n}(X)[0] \ar[r]^-{\sim} & \R \sigma_{n*} \L\sigma_n^* \Z_{\leq n}(X)[0]}$$
is invertible. This follows immediately from $\L\sigma_n^*\Z_{\leq n}(X)[0]=\Z_{tr}(X)[0]$ as $\Z_{\leq n}(X)[0]$ is $\sigma_n^*$-admissible.
\end{proof}

\begin{cor}
\label{derived-strong-ngen}
The functor $\L\sigma_n^*: D(Shv_{tr }^{\tau}(k_{\leq n}))\longby{ }D(Shv_{tr}^{\tau}(k))$ is a fully faithful embedding. It induces an equivalence of triangulated categories between 
$D(Shv_{tr }^{\tau}(k_{\leq n}))$ 
and the triangulated subcategory of $D(Shv_{tr}^{\tau}(k))$ stable under arbitrary sums and generated by the complexes $\Z_{tr}(X)[0]$ for $X\in (Sm/k)_{\leq n}$.
\end{cor}

\begin{proof}
Follows directly from Lemma \ref{Lsigma-up-star}.
\end{proof}

\subsubsection{Motivic complexes} Let $\catM$ be a model category (satisfying some technical assumptions such as being cellular and proper on the left) and $S$ be a set of arrows in $\catM$. Then the Bousfield localization $\L_S(\catM)$ exists. As abstract categories, $\L_S(\catM)=\catM$, the cofibration are the same and $S$ is contained in the class of weak equivalences of $\L_S(\catM)$. Moreover, the identity functor $\catM\to \L_S(\catM)$ is a Quillen functor. This means that $Ho(\catM)\to Ho(\L_S(\catM))$ admits a right adjoint which identifies $Ho(\L_S(\catM))$ with the full subcategory of $Ho(\catM)$ consisting of $S$-local objects (\cf \cite[Th. 4.3.1]{HIR}). 
In other words, we can define $Ho(L_S(\catM))$ (up to an equivalence of categories) as being the full subcategory of $S$-local objects in $Ho(\catM)$. Up to this equivalence of categories, $\L_S : Ho(\catM)\to Ho(\L_S(\catM))$ becomes the localisation functor and is the left adjoint to the inclusion. 

The triangulated category $\DM^{\tau}_{\eff}(k)$ is the homotopy category of a Bousfield localization $\L_S(\catM)$ where $\catM$ is the category of complexes of $\tau$-sheaves with transfers and 
$$S=\{\text{\rm maps of the form $\Z_{tr}(\Aff^1_X)\to \Z_{tr}(X)$ and their shifts}\}$$

Therefore $\DM^{\tau}_{\eff}(k)$ is the full subcategory of $D(Shv^{\tau}_{tr}(k))$ whose objects are the $\Aff^1$-local complexes (called also motivic complexes), \ie these are complexes $\cA_{\d}$ such that:
$$\Hom_{D(Shv^{\tau}_{tr}(k))}(\Z_{tr}(X), \cA_{\d}[m]) \simeq \Hom_{D(Shv^{\tau}_{tr}(k))}(\Z_{tr}(\Aff^1_X), \cA_{\d}[m])$$
We denote by $\L_{\Aff^1}:\xymatrix{D(Shv^{\tau}_{tr}(k)) \ar[r] & \DM^{\tau}_{\eff}(k)}$ the $\Aff^1$-localization functor which is left adjoint to the obvious inclusion. 

For bounded above complexes one can also use \cite[Lect. 14]{VL}.  One can easily see that, for $\tau = $ Nisnevich and $k$ perfect, the resulting triangulated category of bounded above (effective) motivic complexes is fully embedded in $\DM^{\tau}_{\eff}(k)$. In fact, one can use the description of these categories as full subcategories of $D(Shv^{\tau}_{tr}(k))$ (recall that $D^-(Shv^{\tau}_{tr}(k))\subset D(Shv^{\tau}_{tr}(k))$) and just check that a bounded above $\Aff^1$-local object is also an $\Aff^1$-local object of $\DM^{\tau}_{\eff}(k)$ (this is equivalent to say that the complex is bounded and the homology sheaves are homotopy invariants by Voevodsky's theorem on the $\Aff^1$-invariance of cohomology).

The object $\L_{\Aff^1}(\Z_{tr}(X))$ will be denoted by $\M(X)$ for any smooth $k$-variety $X$. This is the homological motive of $X$.

\begin{remark}
\label{rem-DM-sums}
The category $\DM^{\tau}_{\eff}(k)$ admits arbitrary sums. 
Moreover, as a triangulated category with arbitrary sums,
$\DM^{\tau}_{\eff}(k)$ is  generated by $\L_{\Aff^1}(X)$ with $X\in
Sm/k$. 
If $\tau\in\{ {\rm co,Nis}\}$ or $k$ has finite cohomological dimension, the
inclusion $\DM^{\tau}_{\eff}(k) \subset D(Shv_{tr}^{\tau}(k))$ 
commutes with arbitrary sums. This follows easily from the commutation of
${\rm R} \ihom(\Aff^1_k,-)$ with arbitrary sums. 
Moreover, the generators 
$\L_{\Aff^1}(X)$
are compact so that
$\DM^{\tau}_{\eff}(k)$ is compactly generated.
\end{remark}

\begin{defn}
We denote by $\DM_{\leq n}^{\tau}(k)$ the triangulated subcategory of
$\DM_{\eff}^{\tau}(k)$ stable under arbitrary sums generated by $\M(X)$
for $X\in (Sm/k)_{\leq n}$. This is the triangulated category of $n$-motives.
\end{defn}

\begin{remark} 
The functor $\L_{\Aff^1}\circ \L\sigma^*_n:D(Shv^{\tau}_{tr}(k_{\leq n})) \to \DM_{\eff}^{\tau}(k)$
takes values in the subcategory $\DM_{\leq n}^{\tau}(k) \subset \DM_{\eff}^{\tau}(k)$.
\end{remark}

\subsection{The functor $\L \pi_0$} 
\begin{lemma} \label{0-der-equi}
The  functor  $\L\sigma^*_0$ induces an equivalence of triangulated categories $D(Shv^{\tau}_{tr}(k_{\leq 0}))[1/p^{\tau}]\simeq \DM_{\leq 0}^{\tau}(k)[1/p^{\tau}]$ where
$p^{\tau}$ is $1$ unless $\tau=\text{\emph{{\'e}t}}$; in this case, it is the exponential characteristic of the field $k$.
\end{lemma}

\begin{proof} By Corollary \ref{derived-strong-ngen}, the functor
  $\sigma_0^*=\L\sigma_0^*: D(Shv_{tr }^{\tau}(k_{\leq 0}))\longby{
  }D(Shv_{tr}^{\tau}(k))$ is a fully faithful embedding and induces
  an equivalence with the triangulated subcategory 
  $D(Shv_{tr}^{\tau}(k))_{\leq 0}$ of $D(Shv_{tr}^{\tau}(k))$ with arbitrary sums and generated by $\Z_{tr}(l/k)[0]$ with $l$ a finite separable extension of $k$. 

We need only to prove that $D(Shv_{tr}^{\tau}(k))_{\leq 0}$ coincides with $\DM_{\leq 0}^{\tau}(k)$. It is sufficient to show that the objects of $D(Shv_{tr}^{\tau}(k))_{\leq 0}$ are $\Aff^1$-local. 
For this, we remark that any complex $\cA_{\d}$ in 
 $D(Shv_{tr}^{\tau}(k))_{\leq 0}$ 
is the homotopy limits of the bounded complexes 
$\tau_{\leq n} \sigma_{\geq -n} \cA$.
As the property of being $\Aff^1$-local is stable under homotopy limits
we may assume $\cA_{\d}$ to be a bounded complex of $0$-generated sheaves. In fact, we may assume that $\cA_{\d}$ is  concentrated in degree $0$ with value the $0$-motivic sheaf $\cL$.

We are left to show that $\cL$ is strictly $\Aff^1$-invariant.  For $\tau \neq \text{\'et}$ there is nothing to prove as the higher cohomology groups with values in $\cL$ vanish.
For $\tau=\text{\'et}$ the result follows from Proposition \ref{prop-abelian-hitr}.
\end{proof}

\begin{propose}
\label{exist-l-pi-0}
Assume one of these two conditions:
\begin{enumerate}

\item[(a)] $\tau\neq \text{\emph{{\'e}t}}$,

\item[(b)] we work with rational coefficients.

\end{enumerate}
The functor $\pu_0^*$ admits a total left derived functor: 
$$\L\pu_0^*:\xymatrix{D(Shv_{tr}^{\tau}(k)) \ar[r] & D(Shv_{tr}^{\tau}(k_{\leq 0}))}$$
which is left adjoint to 
$\sigma_0^*:\xymatrix{D(Shv_{tr}^{\tau}(k_{\leq 0})) \ar[r] & D(Shv_{tr}^{\tau}(k))}$.
\end{propose}

\begin{proof}
Using Proposition \ref{existance-LF} we need to show the existence of enough
$\pu_0^*$-admissible complexes. We shall prove that $\Z_{tr}(X)$ is
$\pu_0^*$-admissible for any smooth $k$-variety $X$. 
We remark that under one of the above two conditions, $\Z_{tr}(X)$ is
compact. If follows from Lemma \ref{lemma-reduction-inj} that we need only
to check the vanishing of $\Ext^i(\Z_{tr}(X),\pu_{0*}\cI)=0$ for $i>0$ and $\cI$ injective.

The result follows from the vanishing of higher cohomology in any strongly $0$-generated $\tau$-sheaf $\cL$: for $\tau\neq \text{{\'e}t}$ this is clear; for $\tau=\text{{\'e}t}$, {\'e}tale cohomology with value in the $\Q$-sheaf $\cL$ is also zero in higher degrees. 
\end{proof}

\begin{cor}
Under the conditions of Proposition \ref{exist-l-pi-0} the inclusion $\DM_{\leq 0}^{\tau}(k) \subset \DM_{\eff}^{\tau}(k)$
admits a left adjoint 
$$\L\pi_0:\xymatrix{\DM_{\eff}^{\tau}(k) \ar[r] & \DM_{\leq 0}^{\tau}(k)}$$
\end{cor}

\begin{propose}
Under the conditions of Proposition \ref{exist-l-pi-0},
the functor $\L\pi_0$ takes compact objects to compact objects.

\end{propose}

\begin{proof}
This follows formally from the fact that the functor admits a right adjoint that commutes with arbitrary sums.
\end{proof}

\subsection{The functor $\L\Alb$}
In this section we construct the functor $\LAlb$ for non necessarily
constructible (\ie compact or geometric) motives. 
This extends the functor $\LAlb$ constructed in \cite{BK}. In this sub-section we assume that one of the following conditions is fulfilled:
\begin{itemize}

\item $\tau=\text{Nis}$ and the exponential characteristic $p$ of $k$ is inverted,

\item $\tau=\text{\'et}$ and we work with rational coefficients.

\end{itemize}

\begin{thm}
\label{thm-LAlb-main}
Under one of the above conditions, we have:

{\rm (i)} The composition $\iota: \HI_{\leq
  1}^{\et}(k) \subset Shv_{tr}^{\et}(k) \subset Shv_{tr}^{\tau}(k)$ 
admits a right derived
  functor $\R\iota: D(\HI^{\et}_{\leq 1}(k))\subset
  D(Shv^{\tau}_{tr}(k))$.
Moreover, with rational coefficients $\R\iota$ is a full-embedding with
  essential image the subcategory $\DM^{\et}_{\leq 1}(k)$.

{\rm (ii)} The composition: 
$$\xymatrix{Shv^{\tau}_{tr}(k)
  \ar[r]^-{a_{\et}} & Shv^{\et}_{tr}(k) \ar[r]^-{\Alb} &
  \HI^{\et}_{\leq 1}}$$ 
which we still denote by $\Alb$ and which is left adjoint to $\iota$,
admits a total left derived functor
  $\LAlb:\xymatrix{D(Shv^{\tau}_{tr}(k)) \ar[r] & D(\HI_{\leq 1}^{\et})}$
  which is left adjoint to $\R \iota$. 
Moreover, with rational coefficients, the counit of the adjunction $\LAlb \circ \R\iota \by{\sim} id$ is invertible.

{\rm (iii)} The functor $\LAlb$ factors through the $\Aff^1$-localization functor:
$$\xymatrix{D(Shv^{\tau}_{tr}(k)) \ar[r]^-{\LAlb}\ar[d]_-{\L_{\Aff^1}}  & D(\HI^{\et}_{\leq 1}(k)) \\
\DM_{\eff}^{\tau}(k) \ar@{.>}[ur] }$$
The dotted functor will be also called $\LAlb$.
\end{thm}

The last assertion of (ii) implies the
last assertion of (i). 
The existence of the right derived functor $\R\iota$ is clear as 
$\HI_{\leq 1}^{\et}(k)$ is a Grothendieck abelian category and hence admits enough fibrant complexes. 
  
To prove the existence of $\L\Alb$ we will prove that there are enough $\Alb$-admissible complexes in $Shv^{\tau}_{tr}(k)$. We use Fulton's definition of algebraic equivalence and denote $\NS^r(X)$ the group of codimension $r$ cycles modulo algebraic equivalence.

\begin{defn}  
A smooth $k$-scheme $X$ is said to be \emph{$\NS^r$-local} if $\NS^r(X_{k_{s}})=0$ where $k_s/k$ is a separable closure of $k$ and $X_{k_s}=X\otimes_k k_s$.
\end{defn}

\begin{remark}\label{fulton-alg-curve}
When $k$ is separably closed and the exponent characteristic of $k$ is inverted, one can show that $\alpha\in \CH^r(X)$ is algebraically equivalent to zero if and only if there exist a smooth projective curve $C$, two rational points $x_0, \, x_1\in C(k)$ and
$\beta\in \CH^r(C\times_k X)$ such that $(x_0\times {\rm id}_X)^*\beta=0$ and $(x_1\times {\rm id}_X)^*\beta=\alpha$.  
\end{remark}

\begin{propose}
\label{prop-enough-admissible}
Let $X$ be a smooth $k$-scheme which is affine and $\NS^1$-local.  Then $\Z_{tr}(X)$ is $\Alb$-admissible. 
\end{propose}

\begin{proof}
The object $\Z_{tr}(X)$ is compact if $\tau=\text{Nis}$ or if we work
with rational coefficients. By Lemma \ref{lemma-reduction-inj}
we need to check that $\text{H}_{\text{Nis}}^*(X,\cI)=0$ for $*>0$ and
$\cI$ injective in $\HI_{\leq 1}^{\et}(k)$. 
Let $\cL=\ihom(\G_m,\cI)=\cI_{-1}$ be the Voevodsky contraction of $\cI$; by Corollary \ref{cor-V-contract-1-0} this
is a torsion free $0$-motivic 
{\'e}tale sheaf.
Form the exact sequence in $\HI^{\rm Nis}_{tr}(k)$:
$$0\to \cN \to \cL\otimes \G_m \to \cI \to \cK \to 0$$

As $\ihom(\G_m,-)$ is an exact functor, it follows that $\cN$ and $\cK$ are
birational homotopy invariant sheaves with transfers. We deduce that
$\text{H}_{\text{Nis}}^*(X,\cN)=\text{H}_{\text{Nis}}^*(X,\cK)=0$ for
$*>0$. We have also $\text{H}_{\text{Nis}}^*(X,\cL\otimes \G_m)=0$ for $*>1$.
It follows that for $*>1$ one has
$\text{H}_{\text{Nis}}^*(X,\cI)=0$ and we get a surjection:
$$\xymatrix{\text{H}_{\text{Nis}}^1(X,\cL\otimes \G_m) \simeq \text{H}_{\text{Nis}}^1(X,\cL\otimes \G_m/\cN) \ar@{->>}[r] & \text{H}_{\text{Nis}}^1(X,\cI) }$$
Using the Leray spectral sequence $\text{H}^p_{\text{Nis}}(X,\R^q\theta_*
\cF)\Rightarrow \text{H}_{\et}^{p+q}(X,\cF)$ for the morphism of sites
$\theta: X_{\et} \to X_{\text{Nis}}$ we deduce as usual an inclusion
$\text{H}_{\text{Nis}}^1(X,\cI) \subset  \text{H}_{\et}^1(X,\cI)$. In
particular we need only to show that the map:
$$\text{H}_{\et}^1(X,a_{\et}(\cL\otimes \G_m)) \to \text{H}_{\et}^1(X,\cI)$$
is zero.

As $X$ is affine and $\NS^1$-local, by Lemma \ref{lemma-het1=ext1} below one has an isomorphism:
$$\Ext^1_{\HI_{\leq 1}^{\et}(k)} (\Alb(X),a_{\et}(\cL\otimes \G_m)) \simeq \text{H}_{\et}^1(X,a_{\et}(\cL\otimes \G_m)) $$

Consider the commutative square:

$$\xymatrix{ \Ext^1(\Alb(X),a_{\et}(\cL\otimes \G_m)) \ar[r] \ar[d]_-{\sim} & \Ext^1(\Alb(X),\cI)  \ar[d] \\
  \text{H}_{\et}^1(X,a_{\et}(\cL\otimes \G_m)) \ar[r] & \text{H}_{\et}^1(X,\cI) }$$
To conclude, remark that $\Ext^1(\Alb(X),\cI)=0$ since $\Alb(X)$ is a $1$-motivic sheaf and $\cI$ is injective in $\HI^{\et}_{\leq 1}(k)$. 
\end{proof}

\begin{lemma}
\label{lemma-het1=ext1}
Let $X$ be a smooth affine scheme which is $\NS^1$-local. For any $0$-motivic \'etale sheaf $\cL$ which is torsion free, the obvious morphism:
\begin{equation}
\label{eq: Ext1-Het1}
\xymatrix{\Ext^1(\Alb(X),\cL\otimes\G_m) \ar[r] & {\rm H}_{\et}^1(X,\cL\otimes \G_m)}
\end{equation}
is an isomorphism. (Here, we write $\cL\otimes \G_m$ for the tensor product of homotopy invariant \'etale sheaves with transfers, \ie what was written $a_{\et}(\cL\otimes \G_m)$ in the proof of 
Proposition \ref{prop-enough-admissible})
\end{lemma}

\begin{proof}
We break the proof into three steps. The first one is a reduction to the case of a separably closed base field:

\smallskip

\noindent
\emph{Step 1:} Fix $k\subset k_{s}$ a separable closure of $k$ and denote by $G_k$ its Galois group. We assume that the lemma holds over $k_s$, \ie the morphism of $G_k$-modules:  
$$\xymatrix{\Ext^1(\Alb(X_{k_s}),\cL_{k_s}\otimes\G_m) \ar[r] & {\rm H}_{\et}^1(X_{k_s},\cL_{k_s}\otimes \G_m)}$$
is invertible. On the other hand, the universality of the Albanese scheme gives the following isomorphism of $G_k$-modules:
$$\Ext^0(\Alb(X_{k_s}),\cL_{k_s}\otimes\G_m) \to {\rm H}_{\et}^0(X_{k_s},\cL_{k_s}\otimes \G_m)$$

Using the morphism of the two Hochschild-Serre spectral sequences:
$$\xymatrix{ {\rm H}^p(G_k, \Ext^q(\Alb(X_{k_s}),\cL_{k_s}\otimes \G_m)) \ar@{=>}[r] \ar[d] & \Ext^{p+q}(\Alb(X),\cL\otimes \G_m) \ar[d]\\ 
{\rm H}^p (G_k, {\rm H}_{\et}^q(X_{k_s},\cL_{k_s}\otimes \G_m)) \ar@{=>}[r] & {\rm H}^{p+q}_{\et}(X,\cL\otimes \G_m)}$$
we obtain a morphism of exact sequences:
$$\xymatrix{ 0 \ar[d] & 0 \ar[d]\\
\Ext^1(\Alb(X),\cL\otimes \G_m)  \ar[r] \ar[d]&   {\rm H}^1_{\et}(X,\cL\otimes \G_m)  \ar[d] \\
{\rm H}^0(G_k,\Ext^1(\Alb(X_{k_s}), \cL_{k_s}\otimes \G_m))  \ar[r]^-{\sim} \ar[d]&   {\rm H}^0(G_k,{\rm H}^1_{\et}(X_{k_s},\cL_{k_s}\otimes \G_m)) \ar[d] \\
{\rm H}^1(G_k,\Ext^0(\Alb(X_{k_s}), \cL_{k_s}\otimes \G_m))  \ar[r]^-{\sim} &  {\rm H}^1(G_k,{\rm H}^0_{\et}(X_{k_s},\cL_{k_s}\otimes \G_m)) }$$
By the five lemma we are then done.

\smallskip

\noindent
\emph{Step 2:} From now on, we assume our base field $k$ to be separably closed.
$\cL$, being torsion free, is a filtered colimit of free lattices. We may thus assume $\cL$ to be the constant sheaf $\Z$.

First prove the surjectivity of \eqref{eq: Ext1-Het1}, \ie
$$\xymatrix{\Ext^1(\Alb(X),\G_m) \ar[r] & {\rm H}_{\et}^1(X,\G_m)=\Pic(X)}$$
Let $\cE_1$ be a line
bundle on $X$. As $X$ is $\NS^1$-local, we know that the 
class $[\cE_1]\in \Pic(X)$ is algebraically
equivalent to zero. By Remark \ref{fulton-alg-curve} there exist a smooth projective 
curve $C$ with two points $x_0, x_1\in C(k)$ 
and a line bundle $\cE$ on $X\times_k C$ such that $\cE_{|X\times
  x_0}$ is free and $\cE_{|X\times x_1}\simeq \cE_1$.

Let us choose  a trivialization $t:\mathcal{O}_{X\times x_0} \simeq \cE_{|X\times
  x_0}$. We get then an element $(\cE,t)\in \Pic(X\times C, X\times
  x_0)$ which by Voevodsky \cite{Voe1} gives a correspondence
  (unique up to homotopy)   $\alpha\in Cor(X, C-x_0)$. Recall the
  construction of $\alpha$. As $X$ is affine, $X\times x_0$ admits an
  affine neighborhood in $X\times C$ (for example $X\times (C-x)$ for
  any closed point $x\in C$ different from $x_0$). It follows that it
  is possible to extend the
  trivialization $t$ to a trivialization $t':\mathcal{O} \simeq \cE$
  over an open neighborhood of $X\times x_0$. The Cartier divisor $\alpha$
  defined by $t'$ has support inside $X\times (C-{x_0})$. As it is
  closed in $X\times C$, it is proper and affine over $X$. This means
  that $\alpha$ is a finite correspondence from $X$ to $C-x_0$. 

It follows from the construction of $\alpha$ that the image of
$[x_1]\in \Pic(C-x_0)$ along the map $\alpha^*:\Pic(C-x_0) \to
\Pic(X)$
is equal to $[\cE_1]$.

Now, $\alpha$ induces a section $\alpha\in \Alb(C-x_0)(X)$ which by
the universality of the Albanese scheme factors:
$$\xymatrix{X \ar[r] \ar@/^1.7pc/[rr]^{\alpha} & \Alb(X) \ar[r] &
  \Alb(C-x_0)}$$
It is clear that $[\cE_1]$ is the image by:
$$\xymatrix{\Ext^1(\Alb(C-x_0),\G_m) \ar[r] & \Ext^1(\Alb(X),\G_m)
  \ar[r] & \Pic(X) }$$
of the class of the extension:
$$\xymatrix{0 \ar[r] & \G_m \ar[r] & \Alb(C-\{x_0,x_1\}) \ar[r] & \Alb(C-x_0)  \ar[r] & 0 }$$
This proves that $[\cE_1]$ is in the image of $\Ext^1(\Alb(X),\G_m)\to \Pic(X)=\text{H}_{\et}^1(X,\G_m)$.

\smallskip

\noindent
\emph{Step 3:} Finally, we prove the injectivity of 
\eqref{eq: Ext1-Het1} (still for $\cL=\Z$). Suppose given an extension:
$$\xymatrix{0 \ar[r] & \G_m \ar[r] & \mathcal{E} \ar[r] & \Alb(X) \ar[r] & 0}$$
$\mathcal{E}$ is then represented by a commutative group scheme. Suppose that the class of this extension goes to zero by 
\eqref{eq: Ext1-Het1}. This means that
the $\G_m$ torsor $X\times_{\Alb(X)}\mathcal{E}$ splits. Fix a splitting $X\to X\times_{\Alb(X)}\mathcal{E}$ and consider the composition:
$$X\to X\times_{\Alb(X)}\mathcal{E} \to \mathcal{E} $$
By the universality of the Albanese scheme we get a morphism of group schemes $\Alb(X) \to \mathcal{E}$ which is clearly a splitting of our initial extension.
\end{proof}

\begin{cor}
\label{existence-enough-Alb-adm}
$C(Shv_{tr}(k))$ admits enough $\Alb$-admissible objects.
\end{cor}

\begin{proof}
It is sufficient to show that any $k$-variety admits a Zariski
hyper-cover by $\NS^1$-local affine varieties.
As $\xymatrix{\NS^1(U_{k_s} ) \ar@{->>}[r] & \NS^1(V_{k_s})}$ is surjective for any open subscheme $V$ of a smooth $k$-scheme $U$,
it is sufficient to prove that every smooth $k$-variety $X$ can be covered by $\NS^1$-local varieties. 
Choose a system of generators $a_1, \dots, a_n$ of the finitely generated module $\NS^1(X_{k_s})$ with $a_i$ representable by a very ample line bundle $\cL_i$ on $X_{k_s}$. For any point $x\in X$, one can find divisors $D_i\subset X_{k_s}$ representing $\cL_i$ and which are disjoint from $x\otimes_k k_s$. Denote by $D'_i$ the image of $D_i$ by $X_{k_s} \to X$. It follows that $X-\cup_i D'_i$ is an $\NS^1$-local neighborhood of $x$.
\end{proof}

{\it Proof of Theorem \ref{thm-LAlb-main}}.
Corollary \ref{existence-enough-Alb-adm} shows the existence of a left adjoint $\LAlb$ to $\R\iota$ by the general Proposition \ref{existance-LF}. Let us shows that 
$\LAlb$ factors through the $\Aff^1$-localization functor $\L_{\Aff^1}$. For this recall that:
$$\L_{\Aff^1}:\xymatrix{D(Shv^{\tau}_{tr}(k)) \ar[r] & \DM^{\tau}_{\eff}(k)}$$
identify $\DM^{\tau}_{\eff}$ with the Verdier localization of $D(Shv^{\tau}_{tr}(k)) $ with respect to the triangulated subcategory $\cI$ stable by infinite sums and generated by the complexes:
$$\cQ_X=[\xymatrix{0 \ar[r] & \Z_{tr}(\Aff^1_X) \ar[r] & \Z_{tr}(X) \ar[r] & 0  }]$$
Remark that $\cI$ is also generated by $\cQ_X$ with $X$ supposed $\NS^1$-local. Indeed by the proof of Corollary
\ref{existence-enough-Alb-adm}, every smooth $k$-variety admits a Zariski hyper-cover by $\NS^1$-local affine open subvarieties.
By universality it suffices to show that $\LAlb$ sends these complexes to $0$.  The result follows then from the well known fact that $\Alb(\Aff^1_X)=\Alb(X)$.

To finish the proof, we show that the counit $\xymatrix{\LAlb\circ \R\iota \ar[r]^-{\sim }  & id}$ is invertible with rational coefficients. As both $\LAlb$ and $\R\iota$ commutes with arbitrary sums we need only to check that:
$$\xymatrix{\LAlb(h_0(C)) \ar[r]^-{\sim }  & h_0(C)}$$
with $C$ a smooth open curve (use that $\R\iota\simeq \iota$ with rational coefficients). Recall that $\Z_{tr}(C)\to h^{\et}_0(C)$ is an $\Aff^1$-weak equivalence by \cite{V}. 
As every affine smooth curve is $\NS^1$-local we are left to check that $\Alb(C)\simeq h^{\et}_0(C)$, which is clear.
\eproof

\begin{propose}\label{comp-lalab}
With rational coefficients, the functor $\LAlb$ takes compact objects to compact objects.

\end{propose}

\begin{proof}
By the proof of Corollary
\ref{existence-enough-Alb-adm},
every $k$-variety admits a Zariski hypercover by $\NS^1$-local affine open subvarieties. It follows that the triangulated category $\DM_{\eff}^{\rm Nis}(k)_{\Q}$ is compactly generated by  
the motives of affine $\NS^1$-local smooth $k$-schemes $X$. But for such $X$, we have by construction $\LAlb(\M(X))=\Alb(X)$ which is  compact in $D(\HI_{\leq 1}^{\et}(k))$. Indeed, with rational coefficients $\Alb(X)$ is a direct factor of the motive of a smooth curve which is actually
compact in $\DM^{\rm Nis}_{\eff}(k)_{\Q}$. Our claim follows from the fact that the inclusion $D(\HI_{\leq 1}^{\et}(k))\subset \DM^{\rm Nis}_{\eff}(k)_{\Q}$ commutes with infinite sums. 
\end{proof}

\begin{remark}
By Proposition \ref{comp-lalab} 
we have, with rational coefficients, 
a functor $\LAlb:\xymatrix{\DM_{\eff,\gm}^{\et}(k) \ar[r] &
\DM_{\leq 1, \gm}^{\et}(k)}$. This functor coincides with the one defined by a completely different method in \cite[\S 5]{BK}. Indeed, they are both left adjoint to the obvious inclusion.
\end{remark}

\begin{cor}\label{cor-full-embed-1-mot-etal}
Let $i:\HI_{\leq 1}^{\et}(k)\subset Shv_{tr}^{\et}(k)$ be the obvious inclusion. Then $\R i$ is a full embedding (even with $\Z[1/p]$-coefficients).   
\end{cor}

\begin{proof}
With rational coefficients, this follows from Theorem \ref{thm-LAlb-main} as $\R i$ coincides with $\R\iota$ up to the equivalence $\DM_{\eff}^{\et}(k)\simeq \DM_{\eff}^{\rm Nis}(k)$ (still with rational coefficients). By the Suslin rigidity theorem \cite[Th. 7.20]{VL}, the torsion objects of $\HI_{\leq 1}^{\et}(k)$ are simply the $\sigma_0^*$ of torsion \'etale sheaves with transfers on $(Sm/k)_{\leq 0}$.  
It follows from Lemma \ref{0-der-equi} that $\R i$ restricted to torsion objects is a full embedding.
We conclude now using \cite[B.2.4]{BK}.
\end{proof}

\begin{propose}\label{cohdim-leq1}
The cohomological dimension of the $\Z [1/p]$-linear abelian category $\HI_{\leq 1}^{\et}(k)$ is bounded by $2+{\rm cd}(k)$. Moreover, with rational coefficients, this cohomological dimension is $1$.
\end{propose}

\begin{proof}
Let us define ${\rm cd}'(k)$ to be $2+{\rm cd}(k)$ or $1$ if the coefficients ring is $\Z[1/p]$ or $\Q$. 
We need to show that $\Ext^i(\cA,\cB)=0$ 
for $\cA$ and $\cB$ two $1$-motivic sheaves and
$i>{\rm cd}'(k)$. We split the proof into two steps.

\smallskip

\noindent
\emph{Step 1:}
Using the long exact sequences of $\Ext$-groups associated to
$$\xymatrix{0\ar[r] & \cE^0\ar[r] & \cE\ar[r] & \pi_0(\cE)\ar[r] & 0}$$
for $\cE=\cA$ and $\cE=\cB$ we may assume that each of $\cA$ and $\cB$ is either $0$-motivic or connected (we say that a sheaf $\dagger$ is connected if $\pi_0(\dagger)=0$).

The case where $\cA$ and $\cB$ are both $0$-motivic follows immediately from the Hochschild-Serre spectral. We get actually the more precise statement $\Ext^i(\cA,\cB)=0$ for $i>{\rm cd}'(k)-1$.

We now assume that one of the sheaves $\cA$ or $\cB$ is a connected $1$-motivic sheaf. Let $\cE$ be a connected $1$-motivic sheaf and $\cE_{tor}\subset \cE$ its maximal torsion subsheaf. Then by Suslin rigidity theorem \cite[Th. 7.20]{VL} we know that $\cE_{tor}$ is a $0$-motivic sheaf. Moreover, using the fact that $\cE$ is connected, we deduce that 
$\cE'=\cE/\cE_{tor}$ is uniquely divisible (\ie takes values in the category of $\Q$-vector spaces).
Using the long exact sequences of $\Ext$-groups associated to
$$\xymatrix{0 \ar[r] & \cE_{tor} \ar[r] & \cE \ar[r] & \cE' \ar[r] & 0}$$
for $\cE\in \{\cA,\cB\}$ not $0$-motivic, we may assume that 
each of $\cA$ and $\cB$ is either, $0$-motivic or a uniquely divisible connected $1$-motivic sheaf.  
The case where both $\cA$ and $\cB$ are $0$-motivic has just been treated. We may then assume that at least one of $\cA$ or $\cB$ is a uniquely divisible connected $1$-motivic sheaf.

Suppose that $\cA$ is a $0$-motivic sheaf. Then $\cB$ is uniquely $\Q$-divisible and we have $\Ext^i(\cA,\cB)=\Ext^i(\cA\otimes \Q,\cB)$. 
As $\cA$ is a $0$-motivic sheaf, $\cA\otimes \Q$ decomposes of as a direct sum of simple $0$-motivic sheaves of $\Q$-vector spaces $\cA\otimes \Q=\oplus_{\alpha}\cA_{\alpha}$ where $\cA_{\alpha}$ is a direct summand of some $\Q_{tr}(\Spec(l_{\alpha}))$ with $l_{\alpha}/k$ a finite separable extension. 
Using that $\Ext^i(\cA,\cB)=\prod_{\alpha}\Ext^i(\cA_{\alpha},\cB)$ we 
may assume that $\cA=\Q_{tr}(\Spec(l))$ for some finite separable extension $l/k$. But then we get (using Corollary 
\ref{cor-full-embed-1-mot-etal}): 
$$\Ext^i(\Q_{tr}(\Spec(l)),\cB)={\rm H}_{\et}^i(l,\cB)=0$$ 
for $i>0$ (and in particular for $i>{\rm cd}'(k)-1$) as $\cB$ is uniquely divisible. 

\smallskip

\noindent
\emph{Step 2:} By Step 1, we 
may assume that $\cA$ is a uniquely divisible and connected $1$-motivic sheaf. 

Let $\cL=\ihom(\G_m,\cA)$. This is a $0$-motivic sheaf by 
Corollary \ref{cor-V-contract-1-0}.
Consider the exact sequence of \'etale sheaves
$$0\to \cN\to \cL\otimes \G_m\to \cA\to \cA_b\to 0$$
Then $\cN$ is $0$-motivic and $\cA_b$ is a birational, uniquely divisible and connected $1$-motivic sheaf. 
Using the long exact sequence of $\Ext$-groups we 
need to consider the following two cases:

\begin{enumerate}

\item $\cA=\cL\otimes \G_m/\cN$ with $\cL$ and $\cN$ two uniquely divisible $0$-motivic sheaves,

\item $\cA=\cA_b$ is a birational, uniquely divisible and connected $1$-motivic sheaf.

\end{enumerate}
Using other long exact sequences of $\Ext$-groups, one easily sees that (1) and (2) follow from the following properties:

\begin{enumerate}

\item[(i)] If $\cN$ is $0$-motivic and uniquely divisible then $\Ext^i(\cN,\cB)=0$ for $i>{\rm cd}'(k)-1$,

\item[(ii)] If $\cL$ is $0$-motivic and uniquely divisible then $\Ext^i(\cL\otimes \G_m,\cB)=0$ for $i>{\rm cd}'(k)$,

\item[(iii)] If $\cA$ is a birational, uniquely divisible and connected $1$-motivic sheaf then 
$\Ext^i(\cA,\cB)=0$ for $i>{\rm cd}'(k)$.

\end{enumerate}

Property (i) has been proved in Step 1.
For (ii), we can write $\cL$ as a direct sum $\cL=\oplus_{\alpha} \cL_{\alpha}$ where $\cL_{\alpha}$ are direct summand of $\Q_{tr}(\Spec(l_{\alpha}))$ with $l_{\alpha}/k$ finite separable extensions. 
It is then sufficient to show that 
$\Ext^i(\Q_{tr}(\Spec(l))\otimes \G_m,\cB)=0$ for $l/k$ finite and separable and $i>{\rm cd}'(k)$. 
Consider now the exact sequence
$$\xymatrix{0 \ar[r] & \Z_{tr}^{\et}(l)\otimes \mu_{\infty} \ar[r] & \Z_{tr}^{\et}(l)\otimes \G_m \ar[r] & 
\Q_{tr}(l)\otimes \G_m \ar[r] & 0}$$
where we wrote $l$ in place of $\Spec(l)$ and $\mu_{\infty}$ for the torsion subsheaf of $\G_m$. Using the case when $\cA$ is $0$-motivic, settled in Step 1, we are reduced to show that 
$\Ext^i(\Z_{tr}^{\et}(l)\otimes \G_m,\cB)=0$ for $i>{\rm cd}'(k)$. 
Consider now the curve $C_1=(\Aff^1_k-o)\otimes_k l$. The sheaf
$\Z_{tr}^{\et}(l)\otimes \G_m$ is a direct summand of the motive $\M(C_1)$. Using Corollary \ref{cor-full-embed-1-mot-etal}, 
it is sufficient to show that  
$$\Hom_{\DM_{\eff}^{\et}(k)}(\M(C_1),\cB[i])={\rm H}_{\et}^i(C_1,\cB)=0$$
for $i>2+{\rm cd}(k)$ (resp. $i>1$ with rational coefficients). The integral case follows from \cite[X, Cor. 4.3]{SGA4} as $C_1$ has Krull dimension $1$. With rational coefficients, we use that
${\rm H}_{\et}^i(C_1,-)={\rm H}_{\rm Nis}^i(C_1,-)$ and the well known fact that the Nisnevich cohomological dimension is bounded by the Krull dimension. 

For (iii), we use Theorem 
\ref{struct-1-mot} to get an exact sequence
$$\xymatrix{0 \ar[r] & \cT \ar[r] & \cA' \ar[r] & \cA \ar[r] & 0}$$
with $\cT$ a uniquely divisible $0$-motivic sheaf and $\cA'$ a direct sum of abelian varieties tensored by $\Q$. Using the long exact sequence of $\Ext$-groups and the case of $0$-motivic sheaves, settled in Step 1, we may assume that $\cA=\oplus_{\beta} A_{\beta}\otimes \Q$ with $A_{\beta}$ abelian varieties. We then have $\Ext^i(\cA,\cB)=\prod_{\beta} \Ext^i(A_{\beta}\otimes \Q,\cB)$ so we may assume $\cA=A\otimes \Q$ for some abelian variety $A$. One can find an irreducible smooth and projective curve $C_2$ having a rational point $c$ such that $A\otimes \Q$ is a direct factor of 
$h_0^{\et}(C,c)\otimes \Q$. Using the exact sequence
$$\xymatrix{0\ar[r] & h_0^{\et}(C,c)_{tor} \ar[r] & h_0^{\et}(C,c)\ar[r] & h_0^{\et}(C,c)\otimes \Q \ar[r] & 0}$$
and the fact that $h_0^{\et}(C,c)_{tor}$ is a direct factor of $\M(C)$, we reduce to show (by Corollary \ref{cor-full-embed-1-mot-etal}) that 
$$\Hom_{\DM^{\et}_{\eff}(k)}(\M(C_2)\otimes \Q,\cB[i])={\rm H}_{\et}^i(C_2,\cB)=0$$ 
for $i>2+{\rm cd}(k)$ (resp. $i>1$ with rational coefficients). We then argue as for (ii).
\end{proof}

\begin{remark} 
With $\Z[1/p]$-coefficients (and $\tau={\rm Nis}$),  
$\R\iota$ is not the composition of the right
derived functors of the inclusions $i:\HI_{\leq 1}^{\et}(k) \subset
Shv_{tr}^{\et}(k)$ and $j:Shv_{tr}^{\et}(k) \subset Shv_{tr}^{\tau}(k)$.
Let us suppose for simplicity that $k$ is separably closed and pick a prime $\ell$ invertible in $k$. We will prove that 
$\R \iota (\Z/\ell)$ is a bounded complex, whereas $\R j\circ \R i (\Z/\ell)$ is unbounded. 

Let $X$ be an affine smooth and $\NS^1$-local $k$-scheme. We have by adjunction 
${\rm H}^i(\R \Gamma(X,\R\iota (\Z/\ell)))=\Ext^i(\Alb(X),\Z/\ell)$. By Proposition \ref{cohdim-leq1}, these groups vanish for $i>2$. It follows that the complex 
$\R\iota(\Z/\ell)$ is bounded above by $2$ as 
$h_{-i}(\R\iota(\Z/\ell))$ is the Zariski sheaf associated to 
$U\rightsquigarrow {\rm H}^i(\R \Gamma(U,\R\iota (\Z/\ell)))$ and every smooth scheme $U$ can be covered by $\NS^1$-local open affine subschemes.
On the other hand, $\R i(\Z/\ell)=\Z/\ell$ and 
$\R j(\Z/\ell)$ is the object of $\DM_{\eff}^{\rm Nis}(k)$ that represents \'etale cohomology. This object is unbounded. Indeed there are varieties $Y$ of dimension $d$ such that ${\rm colim}_{V\subset Y} {\rm H}^d(\R\Gamma(V,\R j(\Z/\ell))) = {\rm H}^d_{\et}(k(Y),\Z/\ell)\neq 0$. 
\end{remark}

\subsection{The non-existence of a left adjoint for $n\geq 2$}
Here we work with rational coefficients. We take $\tau=\text{Nis}$ and drop the corresponding indexing in the notations. 
A natural generalization of the previous construction is the
following. Consider the smallest triangulated subcategory
$\DM_{\leq n}(k)$ of $\DM_{\eff}(k)$ stable under infinite sums and containing $\M(X)$ for $X$
smooth of dimension $\leq n$. Is there a left adjoint to the obvious
inclusion? Unfortunately, the answer is negative as pointed out (without proof) by
Voevodsky  \cf \cite[\S 3.4]{V}.

In this section we provide a proof of this fact, which is probably
similar to Voevodsky's. Note however that our argument does not
use the motivic conjectural picture.
We assume that such an adjoint exists and denote it by $\Phi_n:\DM_{\eff}(k) \to \DM_{\leq n}(k)$.  We will derive a contradiction. As for the cases $n=0,1$, the functor $\Phi_n$ takes constructible motives to constructible motives. Indeed, the obvious inclusion  $\DM_{\leq n}(k)\subset \DM_{\eff}(k)$ which is right adjoint to $\Phi_n$ commutes with arbitrary sums.
Note the following:

\begin{lemma}\label{lem-non-exist-n}
Assume our base field $k$ is algebraically closed and of infinite transcendence degree over $\Q$. Let $M$ be a constructible motive. If $\Phi_n$ exists then for any finitely generated extension $k\subset K$ the obvious map $\Phi_n(M_K) \to (\Phi_n(M))_{K}$ is invertible.

\end{lemma}

\begin{proof}
Note that the obvious morphism is the one we get by adjunction from the pull-back along $k\subset K$ of $M \to \Phi_n(M)$. 
By replacing $M$ by the cone of $M\to \Phi_n(M)$ we may assume that $\Phi_n(M)=0$. We then need to prove that $\Phi_n(M_K)$=0. 

Consider the universal map $u:M_K \to \Phi_n(M_K)$. As both $M_K$ and $\Phi_n(M_K)$ are constructible, this map is defined over a smooth variety $V$ with generic point $\Spec(K)$. This means that there exists an object $A\in \DM_{\leq n}(V)$ and a morphism $\tilde{u}:M_V \to A$ in $\DM_{\eff}(V)$ whose pull-back to $k(V)$ is  $u$.

Now remark that for any closed point $x\in V$, the pull-back along $x$ of $\tilde{u}$ is a morphism $\tilde{u}_x: M \to A_x$ with $A_x\in \DM_{\leq n}(k)$. As $\Phi_n(M)=0$, the map $\tilde{u}_x$ is necessarily zero. As $k$ has infinite transcendence degree over $\Q$ and because $M$ and $A$ are constructible this implies that $u=0$. This forces $\Phi_n(M_{K})$ to be zero.
\end{proof}

We have:

\begin{cor}\label{cor-non-exist-n}
Assume that $k$ is algebraically closed with infinite transcendence degree over $\Q$.
Let $M$ be a constructible motive. If $\Phi_n$ exists then the obvious morphism $M \to \Phi_n(M)$ induces an isomorphism
$$\ihom(\Phi_n(M), \Z(r)) \to \ihom(M,\Z(r))$$
for $r\leq n$.
\end{cor}

\begin{proof}
To prove this, it suffices to show that for any finitely generated extension $k\subset K$ and any 
$n\in \Z$ the morphism: $$\Hom(\Spec(K),\ihom(\Phi_n(M), \Z(r))) \to \Hom(\Spec(K),\ihom(M,\Z(r)))$$
is invertible. By adjunction and Lemma \ref{lem-non-exist-n}, the above map is the same as:
$$\Hom_{\DM_{\eff}(K)}(\Phi_n(M_K),\Z(r))\to \Hom_{\DM_{\eff}(K)}(M_K,\Z(r))$$
As $\Z(r)$ is in $\DM_{\leq r}(K)\subset \DM_{\leq n}(K)$
this is true by the universality of $M_K\to \Phi_n(M_K)$.
\end{proof}

Having this, it is easy to provide a contradiction. Indeed, for a smooth and projective variety $X$ of dimension $\leq n$ one has 
$$\ihom(\M(X),\Z(n)[2n])\simeq \M(X)(n-\dim(X))[2n-2\dim(X)]$$ 
by \cite[Cor. 4.3.4]{V}.
As the triangulated subcategory of constructible motives in $\DM_{\leq n}(k)$ is generated by motives of smooth and projective varieties of dimension less than $n$, 
we obtain that $\ihom(M,\Z(n))$ is constructible for any constructible object of $\DM_{\leq n}(k)$.
We deduce from Corollary \ref{cor-non-exist-n} that for any constructible motive $M$,  $\ihom(M,\Z(n))$ is constructible.
This is false for $M=\M(X)$ with $X$ a generic quintic in $\mathbb{P}^4$ and $n=2$. 
Indeed, the complex $\ihom(\M(X),\Z(2)[4])$ is concentrated in (homological) positive degree  and its zero homology sheaf $h_0(\ihom(\M(X),\Z(2)[4]))$ is $\text{CH}^2_{/X}$ (see \ref{defn-ch-p-x}). By Theorem \ref{equality-NS-pi-CH}, we get that $\text{L}_0\pi_0(\ihom(\M(X),\Z(2)[4]))\simeq \text{L}_0\pi_0(\text{CH}^2_{/X})=\text{NS}^2_{/X}$. The latter is not finitely generated. For more details, see \cite{Ay-Hub}.

Despite the above negative result, we expect that the following conjecture is true but also quite difficult.

\begin{conj}
With rational coefficients, $\DM_{\leq n}(k)$ is exactly the full subcategory of motivic complexes whose homology sheaves are $n$-motivic in each degree. In other words, the homotopy $t$-structure on $\DM_{\eff}(k)$ restricts to a homotopy $t$-structure on $\DM_{\leq n}(k)$ whose heart is  $\HI_{\leq n}(k)$. Moreover, $\DM_{\leq n}(k)$ has cohomological dimension $\leq n$ with respect to the homotopy $t$-structure, \ie for $\cF$ and $\cG$ $n$-motivic sheaves, we have
$\Hom_{\DM}(\cF,\cG[i])=0$ for $i>n$.
\end{conj}

\section{Computations and applications}

One of the main tasks of this work is to extend the functor $\LAlb$
defined in \cite{BK} to not necessarily constructible motives in order
to apply it to motives of the form $\underline{\Hom}(\Z(n),\M(X))$. Note that such motives are not constructible in general
(\eg $X$ a generic quintic in $\mathbb{P}^4$ and $n=1$, \cf \cite{Ay-Hub}). In 
this section we
use the functors $\text{L} \pi_0$ and
$\LAlb$ to produce some invariants of algebraic varieties. We  begin with some computations.

\subsection{Chow and N\'eron-Severi sheaves} \label{chow-neron-severi}
Let $X$ be a smooth scheme. Recall
that $\CH^r(X)$ denotes the group of codimension $r$ cycles in $X$ up
to rational equivalence. A cycle $\alpha \in \CH^r(X)$ is said
algebraically equivalent to zero if there exist a smooth connected variety $U$,
a zero cycle $\sum_i n_i[u_i]$ in $U$ of degree zero and an element
$\beta \in \CH^r(U\times_k X)$ such that $\alpha=\sum_i n_i (u_i\times {\rm id}_X )_* (u_i\times {\rm id}_X)^*\beta$.

\subsubsection{} \label{defn-ch-p-x}
Recall that the N\'eron-Severi group $\NS^r(X)$ of
codimension $r$ cycles in $X$ is the quotient
$\CH^r(X)/\CH^r(X)_{alg}$ with $\CH^r(X)_{alg}$ the subgroup of
algebraically equivalent to zero cycles. We denote by $\CH^{r,\tau}_{/X}$
the $\tau$-sheaf associated the presheaf
$U\rightsquigarrow \CH^r(U\times_k X)$. 
We define also a $\tau$-sheaf $\NS^{r,\tau}_{/X}$ in the following
way:

\begin{defn}
Suppose that $U$ is connected. 
A cycle $\alpha\in \CH^r(U\times_k X)$ is \emph{algebraically equivalent to
zero relatively to $U$} (or \emph{$U$-algebraically equivalent to zero} for simplicity) if there exist a smooth connected $U$-scheme $V\to U$,
a finite correspondence $\sum_i n_i[T_i] \in Cor(V/U)$ of
degree zero and $\beta \in \CH^r(V\times_k X)$ such that 
$\alpha= \sum_i n_i (t_i\times {\rm id}_X)_* (t_i\times {\rm id}_X)^*\beta$ 
with $t_i$ the finite surjective projection
$T_i\to U$. When $U$ is not connected, we say that $\alpha$ is algebraically equivalent to zero relatively to $U$ if this is the case of the restrictions to $U_0\times_k X$ for $U_0$ any connected component of $U$. 

We denote by $\NS_{/X}^{r,\tau}$ the $\tau$-sheaf associated to the
presheaf 
$$U \rightsquigarrow \CH^r(U\times_k X)/\CH^r(U\times_k
X)_{U-alg}$$ 
where $\CH^r(U\times_k
X)_{U-alg}$ is the subgroup of cycles that are algebraically
equivalent to zero relatively to $U$.
\end{defn}

\begin{propose}
The morphism $\CH^{r,\tau}_{/X} \to \pi_0(\CH^{r,\tau}_{/X})$ factors
uniquely:
$$\xymatrix{\CH^{r,\tau}_{/X} \ar[r] \ar@/^1.7pc/[rr]^-{\;} & \NS^{r,\tau}_{/X} \ar@{.>}[r]^-{s} &
 \pi_0(\CH^{r,\tau}_{/X}) }$$
\end{propose}

\begin{proof}
The uniqueness of $s$ is clear as the first map is surjective. Let us
prove the existence. For this, we need to show that for any smooth $U$
the subgroup $\CH^r(U\times_k X)_{U-alg}$ goes to zero by the map
$\CH^{r,\tau}_{/X}(U) \to \pi_0(\CH^{r,\tau}_{/X}(U))$. Let $\alpha \in
\CH^r(U\times_k X)_{U-alg}$. By definition, there exists a smooth connected
$U$-scheme $V$, a $0$-correspondence $\sum_i n_i[T_i]\in Cor(V/U)$ of degree zero and
an element $\beta \in \CH^r(V\times_k X)$ such that $\beta = \sum_i n_i
(t_i\times {\rm id}_X)_* (t_i\times {\rm id}_X)^* \beta$ with $t_i$ 
the finite surjective projection $T_i\to U$. 

The cycles $\alpha \in \CH^r(U\times_k X)$ and $\beta \in \CH^r(V\times_k X)$
induce morphisms of $\tau$-sheaves:
$$h_0^{\tau}(U) \to \CH^{r,\tau}_{/X} \qquad \text{and} \qquad h_0^{\tau}(V)
\to \CH^{r,\tau}_{/X} $$
Moreover, the finite correspondence $\sum_i n_i[T_i]$ gives a
morphism $h_0^{\tau}(U)\to h_0^{\tau}(V)$ and the equality
 $\alpha = \sum_i n_i
(t_i\times {\rm id}_X)_*(t_i\times {\rm id}_X)^*\beta $ is exactly the
commutativity of the triangle:
$$\xymatrix{h_0^{\tau} (U) \ar[r] \ar[d] &  \CH^{r,\tau}_{/X} \\
h_0^{\tau}(V) \ar[ur] & }$$
We need to show that the composition 
$$\xymatrix{h_0^{\tau} (U) \ar[r]  &  \CH^{r,\tau}_{/X} \ar[r] & \pi_0
(\CH^{r,\tau}_{/X} )}$$
is zero.
For this, we can show that the composition
$$\xymatrix{h_0^{\tau}(U) \ar[r] & h_0^{\tau}(V) \ar[r] & h_0^{\tau}(\pi_0(V))}$$
is zero. This follows immediately from the fact that $\sum_i n_i[T_i]$
is of degree zero and that $V$ is connected. 
\end{proof}

\begin{thm}
\label{equality-NS-pi-CH}
Under one of the following assumptions:

\begin{enumerate}

\item[(a)] $k$ is separably closed and the exponential characteristic of $k$ is inverted,

\item[(b)] $\tau=\text{\emph{\'et}}$ and the exponential characteristic of $k$ is inverted,

\item[(c)] we work with rational coefficients and $\tau={\rm Nis}$,

\end{enumerate}
the morphism $s:\xymatrix{\NS^{r,\tau}_{/X} \ar[r]^-{\sim} &
  \pi_0(\CH^{r,\tau}_{/X})  }$ is invertible.

\end{thm}

\begin{proof}
Remark that it suffices to show that $\NS^{r,\tau}_{/X}$ is a
$0$-motivic sheaf. Indeed, if this is true we get by universality an
inverse $\pi_0(\CH^{r,\tau}_{/X}) \to \NS^{r,\tau}_{/X}$ from the map
$\CH^{r,\tau}_{/X} \to \NS^{r,\tau}_{/X}$.

To check that $\NS^{r,\tau}_{/X}$ is strongly
$0$-generated we might extend the situation to the separable closure of
$k$ using one of the assumptions. Given a smooth variety $U$ we will
show that $\NS^r(X) \to \CH^r(X\times_k U)/\CH^r(X\times_k U)_{U-alg}$
is an isomorphism. This map is obviously injective as it has a
section given by any rational point of $U$. We will show the
surjectivity, \ie every $[Z]\in \CH^r(U\times_k X)$ is
$U$-algebraically equivalent to a ``constant cycle". For this, fix a
point $u\in U$ and consider $\xymatrix{V=U\times_k U \ar[r]^-{pr_1} &
  U}$ together with the finite correspondence of degree zero
$[\Delta]-[U\times u]\in Cor(V/U)$. If $[W]=pr_2^*[Z]\in \CH^r(U\times
U \times X)$ then we have a $U$-algebraically equivalent to zero
cycle:
$$[W\cap (\Delta\times X)]-[W\cap (u\times U \times X)]=[Z]-[U\times
Z_u]$$
This shows that $[Z]$ is $U$-algebraically equivalent to $[U\times
Z_u]$.
\end{proof}

\subsection{The higher N\'eron-Severi groups} 
Here we work only with the Nisnevich topology. We will write $\DM_{\eff}(k)$ instead of $\DM_{\eff}^{\rm Nis}(k)$.

\begin{defn}
Let $X$ be a smooth
$k$-scheme. We define a family of abelian groups $\NS^r(X,s)$ by:
$$\NS^r(X,s)\df \left\{ \begin{array}{ccc} \L_s\pi_0\ihom(\M(X),\Z(r)[2r])(k)  & {\rm for} & r\geq 0,  \\
\\
0 &{\rm for} &r<0. 
\end{array} \right. $$
\end{defn}

\begin{lemma}
For $r>\dim(X)$ we have $\NS^r(X,s)=0$. Moreover, under one of the following hypotheses:
\begin{enumerate}

\item[(a)] $k$ is algebraically closed,

\item[(b)] $k$ is separably closed and the exponential characteristic of $k$ is inverted,

\item[(c)] we work with rational coefficients,

\end{enumerate}
there is a canonical isomorphism $\NS^r(X,0)\simeq \NS^r(X)$ with the
usual N\'eron-Severi group. 
\end{lemma}

\begin{proof}
To prove the vanishing of $\NS^r(X,s)=0$ for $r>\dim{X}=d$ we remark
that $\ihom(\M (X),\Z(r)[2r])\simeq \ihom(\M(X),\Z(d)[2r])\otimes \Z(r-d)$. So
it suffices to show more generally that $\L\pi_0(M\otimes \Z(1))=0$
for any motive $M$. We are reduced to check this for $M=\M(U)$ with $U$
smooth. The result follows then from the fact that
$\pi_0(U)=\pi_0(U\times (\Aff^1_k-o))$.

The complex $\ihom( \M(X),\Z(r)[2r])$ is concentrated in
positive homological degree, \ie the homology sheaf 
$h_i^{\rm Nis}(\ihom(\M (X),\Z(r)[2r]))=0$ for $i<0$. Moreover, we know that
$h_0^{\rm Nis}(\ihom(\M (X),\Z(r)[2r]))$ is the Nisnevich sheaf $\text{CH}_{/X}^r$ associated to the
presheaf $U\rightsquigarrow \text{CH}^r(U\times X)$. We thus have:
$$\L_0\pi_0\ihom(\M (X),\Z(r)[2r])=\pi_0(\text{CH}_{/X}^r)$$
So we need only to show that 
$\pi_0(\text{CH}_{/X}^r)=\NS^r_{/X}$ which is true by Theorem 
\ref{equality-NS-pi-CH}. 
\end{proof}

\begin{propose}
\label{NS-gysin}
Given a closed embedding of smooth schemes $Y\subset X$ of pure
codimension $c$ we have a long exact sequence:
$$\NS^r(X-Y,s+1) \to  \NS^{r-c}(Y,s) \to   \NS^r(X,s) \to
  \NS^r(X-Y,s) $$
\end{propose}

\begin{proof}
This follows immediately from the exact triangle:
$$\xymatrix{\M(X-Y) \ar[r] &  \M(X) \ar[r] & \M(Y)(c)[2c] \ar[r] & }$$
in $\DM_{\eff}^{\rm Nis}(k)$.
\end{proof}

\begin{lemma}
There exists a morphism
${\rm CH}^r(X,s)\to \NS^r(X,s)$ natural in $X$ and compatible with the long exact sequences of Proposition 
\ref{NS-gysin}. 

\end{lemma}

\begin{proof}
By \cite[Prop. 4.2.9 and Th. 4.3.7]{V}, we have  
$${\rm CH}^r(X,s)=h_s^{\rm Nis}(\ihom(\M(X),\Z(r)[2r]))(k)$$
The morphism of the lemma is obtained by applying $h_s^{\rm Nis}$ to 
$$\ihom(\M(X),\Z(r)[2r])\to \L \pi_0 \ihom(\M(X),\Z(r)[2r])$$ 
and then taking $k$-rational points.
\end{proof}

\begin{defn}
We can also define a homological version:
$$\NS_r(X,s)\df  \left\{ \begin{array}{ccc} \L_s\pi_0\ihom(\Z(r)[2r],\M(X))(k)  & {\rm for} & r\geq 0,\\
\\
\L_s\pi_0 (\M(X)\otimes \Z(-r)[-2r])(k) & {\rm for} & r<0. \end{array} \right.$$
\end{defn}

\begin{remark}\label{BOrmk}
Using the formalism of the Grothendieck six operations (\cf \cite{Ayoub}) we think it is possible to extend $\NS^r(X,s)$ to a cohomology theory with support ${\rm H}_{\NS, Z}^*(X,\cdot )$ together with a Borel-Moore homology theory ${\rm H}_{\NS,*}(X,\cdot )$  and a pairing such that these data satisfy the Bloch-Ogus axioms \cite{Bloch-Ogus}. 
 In particular, we would have a Gersten resolution for $\NS^{r}(X,s)$ (\cf \cite{Bloch-Ogus, colliot-kahn}). This deserves a separate treatment. 
\end{remark}

\subsection{The higher Picard and Albanese $1$-motivic sheaves}
Here we still work with $\DM^{\rm Nis}_{\eff}(k)=\DM_{\eff}(k)$.

\begin{defn}
Let $X$ be a $k$-scheme. Define the higher Picard sheaves by:
$$\Pic^r(X,s)\df \left\{ \begin{array}{ccc} \L_s\Alb \, \ihom(\M(X),\Z(r)[2r])  & {\rm for} & r\geq 0,  \\
\\
0 &{\rm for} &r<0. 
\end{array} \right. $$
These are objects of $\HI_{\leq 1}^{\et}(k)$.
\end{defn}

\begin{propose}
\label{vanishing-higher-pic}
We have $\Pic^r(X,*)=0$ when $r>\dim(X)+1$ or if $r=1$ and $*\neq 0, \, 1$. 
Moreover if $k$ is algebraically closed and $X$ smooth,
then 
$\Pic^1(X,0)(k)$ is the usual Picard group of $X$. If $X$ is also projective then $\Pic^1(X,0)$ is represented by the Picard scheme of $X$.
\end{propose}

\begin{proof}
For $r>\dim(X)+1$ we have 
$$\ihom(\M(X),\Z(r))=\ihom(\M(X),\Z(\dim(X))\otimes \Z(r-\dim(X))$$ 
It is then sufficient to show that for any $M\in \DM_{\eff}(k)$ we have $\LAlb(M\otimes \Z(2))=0$. 

We may assume that $M=\M(U)$ with $U$ affine and $\NS^1$-local. 
The result follows then from the
decomposition:
$$\Alb(U\times (\Aff^1_k-0)\times (\Aff^1_k-0))=\Alb(U) \oplus \Alb(U\times (\Aff^1_k-0,1))$$
$$\oplus  \Alb(U\times (\Aff^1_k-0,1)) \oplus \Alb(U\times (\Aff^1_k-0,1)^{\wedge
  2})$$
and the fact that $\Alb(V\times V')=\Alb(V)\otimes \pi_0(V')\oplus
\pi_0(V)\otimes \Alb(V')$ which in
particular implies that $\Alb(U\times (\Aff^1_k-0,1))=\pi_0(U)\otimes \G_m$.
\end{proof}

\begin{propose}
Given a closed embedding of smooth schemes $Y\subset X$ of pure
codimension $c$ we have a long exact sequence:
$$\Pic^r(X-Y,s+1) \to  \Pic^{r-c}(Y,s) \to   \Pic^r(X,s) \to
  \Pic^r(X-Y,s) $$
in $\HI_{\leq 1}^{\et}(k)$.
\end{propose}

\begin{proof}
Same proof as Proposition \ref{NS-gysin}.
\end{proof}

\begin{defn}
For a smooth scheme $X$ define the higher Albanese sheaves by:
$$\Alb_r(X,s)\df \left\{ \begin{array}{ccc}  \L_s\Alb (\ihom(\Z(r)[2r],\M(X)))
& {\rm for} &  r\geq 0 \\
\\
 \L_s\Alb(\M(X)\otimes \Z(-r)[-2r]) & {\rm for} & r<0. \end{array} \right.$$ 
These are objects of $\HI_{\leq 1}^{\et}(k)$.
\end{defn}

\begin{propose}
We have $\Alb_r(X,s)=0$ for $r<-1$ or $r=-1$ and $s\neq 0$. Moreover $\Alb_0(X,0)$ is the usual Albanese scheme $\Alb(X)$ and $\Alb_{-1}(X,0)=\G_m\otimes \pi_0(X)$.
\end{propose}

\begin{proof}
To prove the vanishing of $\Alb_r(X,s)=0$ for $r<-1$ we argue as for Proposition \ref{vanishing-higher-pic}.
\end{proof}

\begin{remark}
Assume $X$ projective of dimension $d$. Then one has
$\Alb_r(X,0)=\Alb(h_0\ihom(\Z(r)[2r],\M(X)))=\Alb(\CH^{d-r}_{/X})$
which most probably, over $k=\C$, will be providing Walter's morphic Abel-Jacobi
map (\cf \cite{Walter}) on the $r$-dimensional cycles.
\end{remark}

\subsection*{Acknowledgments} The first author would like to thank the University of Padova and the European Community's ``Marie Curie early stage researcher" program (under the $6^{th}$ Framework Program RTN {\it Arithmetic Algebraic Geometry}) for hospitality and financial support. 

Both authors are grateful to the referee and the editor for several suggestions which helped them to improve the exposition.

\end{document}